\RequirePackage{fix-cm}
\documentclass[smallextended]{svjour3}       
\smartqed  
\usepackage{graphicx}
\usepackage{tikz}
\usepackage{bm}

\RequirePackage[OT1]{fontenc}
\RequirePackage{amsmath}

\RequirePackage[colorlinks,citecolor=blue,urlcolor=blue]{hyperref}

\usepackage{graphicx}
\usepackage{amsmath}
\usepackage{graphics}
\usepackage{epsfig}
\usepackage{amssymb}
\usepackage{color}
\usepackage{float}
\usepackage[section]{placeins}
\usepackage{subcaption,graphicx}
\usepackage{bbm}

\newcommand{\be}{\begin{eqnarray}}
\newcommand{\ee}{\end{eqnarray}}
\newcommand{\bea}{\begin{eqnarray*}}
\newcommand{\eea}{\end{eqnarray*}}

\newcommand{\norm}[1]{\left\lVert#1\right\rVert}

\newcommand{\AZ}{\textcolor{black}}


\begin{document}

\title{Efficient quantization and weak covering of high dimensional cubes
}


\author{      Jack Noonan      \and Anatoly Zhigljavsky
}


\institute{J. Noonan \at
               School of Mathematics, Cardiff University, Cardiff, CF244AG, UK \\
              \email{Noonanj1@cardiff.ac.uk}
           \and
       A. Zhigljavsky \at
              School of Mathematics, Cardiff University, Cardiff, CF244AG, UK \\
              \email{ZhigljavskyAA@cardiff.ac.uk}           
}

\date{Received: date / Accepted: date}

\maketitle
\begin{abstract}
Let $\mathbb{Z}_n = \{Z_1, \ldots, Z_n\}$ be  a design; that is, a collection  of $n$ points $Z_j \in [-1,1]^d$.
We study the quality  of quantization of $[-1,1]^d$ by the points of $\mathbb{Z}_n$ and the problem of quality of coverage of $[-1,1]^d$
by ${\cal B}_d(\mathbb{Z}_n,r)$, the union of  balls centred at $Z_j \in \mathbb{Z}_n$. We concentrate on the cases where the dimension  $d$ is not small ($d\geq 5$) and $n$
is not too large, $n \leq 2^d$.
We define the design ${\mathbb{D}_{n,\delta}}$ as
a $2^{d-1}$ design defined on vertices of the cube  $[-\delta,\delta]^d$, $0\leq \delta\leq 1$. For this design,
we derive a closed-form expression for the quantization error and very  accurate approximations for {the coverage area}  vol$([-1,1]^d \cap {\cal B}_d(\mathbb{Z}_n,r)) $.
We provide results of a large-scale numerical investigation confirming the accuracy of the developed approximations and the efficiency of the
 designs ${\mathbb{D}_{n,\delta}}$.

\end{abstract}

\section{Introduction}

\subsection{Main notation}

\begin{itemize}

\item $\| \cdot \|$: the Euclidean norm;
 \item ${\cal B}_d(Z,{ r })= \{ Y \in \mathbb{R}^d: \| Y-Z \| \leq { r } \}$: $d$-dimensional ball of radius $r$ centered at $Z \in \mathbb{R}^d$;
 \item $\mathbb{Z}_n = \{Z_1, \ldots, Z_n\}$: a design; that is, a collection  of $n$ points $Z_j \in \mathbb{R}^d$;

 \item ${\cal B}_d(\mathbb{Z}_n,r)= \bigcup_{j=1}^n {\cal B}_d(Z_j,r)$;
 \item \mbox{$C_d(\mathbb{Z}_n,r)= $vol$([-1,1]^d \cap {\cal B}_d(\mathbb{Z}_n,r))/2^d$}: the proportion of the cube $[-1,1]^d$ covered by ${\cal B}_d(\mathbb{Z}_n,r)$;
      \item vectors in $\mathbb{R}^d$ are row-vectors;
 \item  for any $a \in \mathbb{R}$, $\boldsymbol{a} = (a,a,\ldots,a) \in \mathbb{R}^d$.

\end{itemize}

\subsection{Main problems of interest}
\label{main}

We will study the following two main characteristics of designs $\mathbb{Z}_n=\{Z_1, \ldots, Z_n\} \subset \mathbb{R}^d$.

 \begin{bf}{1. Quantization error.}\end{bf} Let  $X=(x_1, \ldots, x_d)$ be uniform random vector on $[-1,1]^d$.
The mean squared quantization error for a design $\mathbb{Z}_n$ is defined by
\be
\label{eq:errorQ}
\theta(\mathbb{Z}_n)=\mathbb{E}_X \varrho^2(X,\mathbb{Z}_n)\,, \;\;{\rm where}\;\;
\varrho^2(X,\mathbb{Z}_n)= \min_{Z_i \in \mathbb{Z}_n} \|X-Z_i\|^2\, .
\ee

\begin{bf}{ {2. Weak covering.}}\end{bf} Denote the proportion of the cube $[-1,1]^d$ covered by the union of $n$ balls
${\cal B}_d(\mathbb{Z}_n,r)= \bigcup_{j=1}^n {\cal B}_d(Z_j,r)$
 by
\bea
\mbox{$C_d(\mathbb{Z}_n,r):= $vol$\left([-1,1]^d \cap {\cal B}_d(\mathbb{Z}_n,r)\right)/2^d$}\,.
\eea
For  given radius $r>0$, the union of  $n$ balls ${\cal B}_d(\mathbb{Z}_n,r)$ makes the
 $(1-\gamma)$-coverage of the cube $[-1,1]^d$ if
\be
\label{weak}
C_d(\mathbb{Z}_n,r)=1-\gamma \,.
\ee
Complete coverage corresponds to $\gamma=0$.  In this paper,  the complete coverage of $[-1,1]^d$ will not be enforced and we
will mostly be interested in {\it weak covering}, that is, achieving  \eqref{weak} with some small $\gamma>0$.\\

Two $n$-point  designs $\mathbb{Z}_n$ and $\mathbb{Z}_n^{\prime}$ will be differentiated in terms of performance as follows:
(a)  $\mathbb{Z}_n$ dominates   $\mathbb{Z}_n^{\prime}$ for  quantization if  $\theta( \mathbb{Z}_n)<\theta( \mathbb{Z}_n^{\prime})$; (b)
if for a given $\gamma \geq 0$,  $C_d(\mathbb{Z}_n,r_1)=C_d(\mathbb{Z}_n^{\prime},r_2)=1-\gamma$ and $r_1<r_2$, then the design $\mathbb{Z}_n$ provides a more efficient $(1-\gamma)$-coverage than $\mathbb{Z}_n^{\prime}$  and is therefore preferable. In Section~\ref{sec:renorm} we extend these definitions by allowing  the two designs to have different number of points and, moreover, to have different dimensions.

\AZ{Numerical construction of $n$-point designs with moderate values of $n$ with good quantization and coverage properties 
has recently attracted much attention in view of diverse applications in several fields including computer experiments \cite{pronzato2012design,pronzato2017minimax,santner2003design}, global optimization \cite{zhigljavsky2021bayesian}, function approximation \cite{SchabackW2006,Wendland2005}
 and numerical integration \cite{pronzato2020bayesian}. Such designs are often referred to as {\it space-filling designs}. Readers can find many additional references in the citations above. Unlike the exiting literature on space-filling, we concentrate on theoretical properties of a family of very efficient designs and derivation of accurate approximations for the characteristics of interest.}

\subsection{Relation between quantization and weak coverage}
\label{sec:relation}

The two characteristics, $C_d(\mathbb{Z}_n,r)$ and $\theta(\mathbb{Z}_n)$, are related: $C_d(\mathbb{Z}_n,r)$, as a function of $r\geq 0$, is the c.d.f. of the r.v. $\varrho(X,\mathbb{Z}_n)$ while $\theta(\mathbb{Z}_n)$ is the second moment of the distribution with this c.d.f.:
\be\label{eq:relation}
\theta(\mathbb{Z}_n)= \int_{r\geq 0} r^2 d  C_d(\mathbb{Z}_n,r)\, .
\ee

In particular, this yields that if an $n$-point design $ \mathbb{Z}_n^\ast$ maximizes, in the set of all $n$-point designs, $C_d(\mathbb{Z}_n,r)$ for all $r>0$, then it also minimizes
$\theta(\mathbb{Z}_n)$. Moreover, if r.v. $\varrho(X,\mathbb{Z}_n)$ stochastically
dominates $\varrho(X,\mathbb{Z}_n^\prime)$, so that $C_d(\mathbb{Z}_n^\prime,r)\leq C_d(\mathbb{Z}_n,r)$ for all $r \geq 0$ and the inequality is strict for at least one $r$, then $\theta(\mathbb{Z}_n) < \theta(\mathbb{Z}_n^\prime)$.

The relation \eqref{eq:relation} can alternatively be written as
\be\label{eq:relation1}
\theta(\mathbb{Z}_n)= \int_{r\geq 0} r \,d C_d(\mathbb{Z}_n,\sqrt{r})\, ,
\ee
where $  C_d(\mathbb{Z}_n,\sqrt{r})$, considered as a function of $r$, is the c.d.f. of the r.v. $\varrho^2(X,\mathbb{Z}_n)$ and hence $\theta(\mathbb{Z}_n)$ is the mean of this r.v. Relation \eqref{eq:relation1} is simply another form of \eqref{eq:errorQ}.

\subsection{{Renormalised versions and formulation of optimal design problems}}
\label{sec:renorm}
In view of \eqref{eq:norm_Qd}, the naturally defined re-normalized version of $\theta(\mathbb{Z}_n)$ is
$Q_d (\mathbb{Z}_n)= {n^{2/d}}  \theta(\mathbb{Z}_n)/({4d} ) .
$
 From \eqref{eq:relation1} and \eqref{eq:relation}, $Q_d (\mathbb{Z}_n)$ is the expectation of the r.v.
${n^{2/d}} \varrho^2(X,\mathbb{Z}_n)/({4d} )$
and the second moment of the r.v.
${n^{1/d}} \varrho(X,\mathbb{Z}_n)/({ 2 \sqrt{d}} )$ respectively. This suggests the following  re-normalization of the radius $r$ with respect to $n$ and $d$:
\be
\label{eq:Rr}
R ={n^{1/d}} r/({ 2 \sqrt{d}} )\, .
\ee
We can then define optimal designs as follows.
Let $d$ be fixed, ${\cal Z}_n= \{\mathbb{Z}_n  \} $ be the set of all $n$-point designs and ${\cal Z}= \cup_{n=1}^\infty {\cal Z}_n$ be the set of all designs.

\begin{definition}
The design $\mathbb{Z}_m^{\ast}$ with some $m$ is optimal for quantization in $[-1,1]^d$, if
\be \label{opt_des_quant}
Q_d (\mathbb{Z}_m^{\ast}) = \min_{n} \min_{\mathbb{Z}_n \in {\cal Z}_n }  Q_d (\mathbb{Z}_n) =  \min_{\mathbb{Z} \in {\cal Z} }  Q_d (\mathbb{Z})\, .
\ee
\end{definition}

\begin{definition}

 The design $\mathbb{Z}_m^{\ast}$ with some $m$ is optimal for $(1-\gamma)$-coverage of  $[-1,1]^d$, if
\be
\label{eq:thick4}
R_{1\!-\!\gamma}  (\mathbb{Z}_m^{\ast}) = \min_{n} \min_{\mathbb{Z}_n \in {\cal Z}_n } R_{1\!-\!\gamma}  (\mathbb{Z}_n) =  \min_{\mathbb{Z} \in {\cal Z} }  R_{1\!-\!\gamma}  (\mathbb{Z})\, .
\ee
Here $0 \leq \gamma \leq 1$ and  for a given design $\mathbb{Z}_n \in {\cal Z}_n$,
\be \label{eq:R_g}
 R_{1\!-\!\gamma} (\mathbb{Z}_n)={n^{1/d}} r_{1\!-\!\gamma}(\mathbb{Z}_n )/({ 2 \sqrt{d}} )\, ,
 \ee
  where
 $r_{1\!-\!\gamma}(\mathbb{Z}_n )$ is defined as the smallest $r$ such that  $C_d(\mathbb{Z}_n,r)=1-\gamma$.

\end{definition}

Importance of the factor $\sqrt{d}$ in \eqref{eq:Rr} will be seen in Section~\ref{sec:asymptotics} where we shall study  the asymptotical behaviour of $(1-\gamma)$-coverings for large $d$.

\subsection{Thickness of  covering}
Let $\gamma=0$ in Definition 2. Then $r_{1}(\mathbb{Z}_n )$  is the covering radius associated with $\mathbb{Z}_n $ so that the union of the balls
${\cal B}_d(\mathbb{Z}_n,r)$ with $r=r_1(\mathbb{Z}_n )$ makes a coverage of $[-1,1]^d$. Let us tile up the whole space $\mathbb{R}^d$ with
the translations of the cube $[-1,1]^d$ and corresponding translations of the balls ${\cal B}_d(\mathbb{Z}_n,r)$. This would make a full coverage of the whole space; denote this space coverage by  ${\cal B}_d(\mathbb{Z}_{(n)},r)$. The  thickness $\Theta$ of any space covering is defined,
see \cite[f-la (1), Ch. 2]{Conway}, as the average number of balls containing a point
of the whole space. In our case of  ${\cal B}_d(\mathbb{Z}_{(n)},r)$, the thickness is
\bea
\Theta(  {\cal B}_d(\mathbb{Z}_{(n)},r)) =\frac{n \, \mbox{vol}\left({\cal B}_d(0,r)\right)}{\mbox{vol}([-1,1]^d)}= \frac{n \,r^d \,
\mbox{vol}\left({\cal B}_d(0,1)\right)}{2^d}\, .
\eea
The normalised thickness, $\theta$, is the thickness $\Theta$ divided by $\mbox{vol}\left({\cal B}_d(0,1)\right)$, the volume of the unit ball, see \cite[f-la (2), Ch. 2]{Conway}. In the case of ${\cal B}_d(\mathbb{Z}_{(n)},r)$, the normalised thickness is
\bea
\theta( {\cal B}_d(\mathbb{Z}_{(n)},r)) = \frac{n \,r^d }{2^d}=d^{d/2}\, \left[R_1(\mathbb{Z}_{(n)} )\right]^d \,\, ,
\eea
where we have recalled that $r=r_1(\mathbb{Z}_n )$ and $ R_{1\!-\!\gamma} (\mathbb{Z}_n)={n^{1/d}} r_{1\!-\!\gamma}(\mathbb{Z}_n ) /({ 2 \sqrt{d}} )$
for any $0\leq \gamma \leq 1$.

We  can thus define the normalised thickness of the covering of the cube by the same formula and extend it to  any $0\leq \gamma \leq 1$:
\begin{definition}
 Let ${\cal B}_d(\mathbb{Z}_{n},r)$ be a $(1-\gamma)$-coverage of the cube $[-1,1]^d$ with $0\leq \gamma \leq 1$. Its normalised thickness
is defined by
\be
\label{eq:Thick5}
\theta( {\cal B}_d(\mathbb{Z}_{n},r)) =  (\sqrt{d}R)^d\, ,
\ee
where  $R=n^{1/d}r/({ 2 \sqrt{d}} )$, see \eqref{eq:Rr}.
\end{definition}

In view of \eqref{eq:Thick5}, we can reformulate the definition \eqref{eq:thick4} of the $(1-\gamma)$-covering optimal design by saying that this design
minimizes (normalised) thickness in the set of all $(1-\gamma)$-covering  designs.

\subsection{The design of the main interest}
\label{sec:main_interest}

We will be mostly interested in the following $n$-point design $\mathbb{Z}_n ={\mathbb{D}_{n,\delta}}$ defined only for  $n=2^{d-1}$:\\

\noindent
{\bf Design ${\mathbb{D}_{n,\delta}}$: } {\it
a $2^{d-1}$ design defined on vertices of the cube  $[-\delta,\delta]^d$, $0\leq \delta\leq 1$.
}\\

For theoretical comparison with design  ${\mathbb{D}_{n,\delta}}$, we shall consider the following simple design, which
extends to the integer point lattice $Z^d$ (shifted by ${\boldsymbol{\frac12}}$) in the whole space $\mathbb{R}^d$:
\\

\noindent
{\bf Design ${\mathbb{D}_{n}^{(0)}}$: } {\it
the collection of   $2^{d}$ points $(\pm \frac12, \ldots, \pm \frac12)$, all  vertices of the cube  $[-\frac12,\frac12]^d$.
}\\

Without loss of generality, while considering the design ${\mathbb{D}_{n,\delta}}$ we assume that the point $Z_1 \in{\mathbb{D}_{n,\delta}}= \{Z_1, \ldots, Z_n\}$
is $Z_1= {\boldsymbol{\delta}} =(\delta, \ldots, \delta)$. Similarly, the first point in
 ${\mathbb{D}_{n}^{(0)}}$ is $Z_1= {\boldsymbol{\frac12}} =(\frac12, \ldots, \frac12)$.
Note also that for numerical comparisons, in Section~\ref{sec:numerical} we shall introduce one more design.

The design ${\mathbb{D}_{n,{ 1/2}}}$   extends to the lattice $D_d$ (shifted by ${\boldsymbol{\frac12}}$) containing points $X=(x_1,\ldots, x_d)$ with integer components satisfying $x_1+\ldots+ x_d=0\; (\!\!\!\!\mod 2)$, see \cite[Sect. 7.1, Ch. 4]{Conway}; this lattice is sometimes called `checkerboard lattice'.
The motivation to theoretically study the design ${\mathbb{D}}_{n,\delta}$ is a consequence of numerical results reported  in~\cite{us} and \cite{second_paper}, where 
the present authors have considered $n$-point designs in $d$-dimensional cubes
providing good coverage and quantization and have shown that for all dimensions $d\geq 7$, the design ${\mathbb{D}}_{n,\delta}$ with suitable $\delta$ provides the best quantization and coverage per point among all other designs considered.  Aiming at practical applications mentioned in Section~\ref{main}, our aim was to consider the designs with $n$ which is not too large and in any case does not exceed $2^d$.

%

\AZ{If the number of points $n$ in a design is much larger than $2^d$, then we may use the following scheme of construction of efficient  quantizers in the cube $[-1 , 1]^d$: (a)   construct one of the very efficient lattice space quantizers,  see \cite[Sect. 3, Ch. 2]{Conway}, (b) take the lattice points belonging to a very large cube,  and (c) scale the chosen large cube  to $[-1 , 1]^d$. In view of Theorem 8.9 in \cite{graf2007foundations},  as $n \to \infty$, the normalised quantization error $Q_d (\mathbb{Z}_n)$ of the sequence of resulting designs $\mathbb{Z}_n$ converges to the
respective quantization error of the lattice space quantizer. However, for any given $n$ the study of quantization error of such designs is difficult (both, numerically and theoretically) as there could be several non-congruent types of Voronoi cells  due to boundary conditions.
Note also that the boundary conditions make significant difference in relative efficiencies of the resulting designs. In particular, the checkerboard lattice $D_d$ is better than the integer-point lattice  $\mathbb{Z}^d$ for all $d\geq 3$  as a space quantizer and becomes the best lattice space quantizer  for $d= 4$ but in the case of cube $[-1,1]^d$, the design
${\mathbb{D}_{n,\delta}}$ (with optimal $\delta$) makes a better quantizer than ${\mathbb{D}_{n}^{(0)}}$ for $d\geq 7$  only; see Section~\ref{quant_exact_results} for theoretical and numerical comparison of the two designs.
}

\subsection{Structure of the rest of the paper and the main results}

In Section~\ref{sec:quant} we study $Q_d (\mathbb{D}_{n,\delta})$, the normalized mean squared quantization error for the design ${\mathbb{D}_{n,\delta}}$. There are two important results, Theorems~\ref{th:1} and~\ref{quant_theorem}. In Theorem~\ref{th:1}, we derive the explicit form for the Voronoi cells for the points of the design ${\mathbb{D}_{n,\delta}}$ and in Theorem~\ref{quant_theorem} we derive a closed-form expression for $Q_d (\mathbb{D}_{n,\delta})$ for any $\delta>0$. As a consequence, in Corollary~\ref{cor:1} we  determine the optimal value of $\delta$.

The main result of Section~\ref{Main_results} is Theorem~\ref{Main_theorem}, where we derive closed-form expressions (in terms of $C_{d,Z,{ r }}$, the fraction of the cube $[-1,1]^d$ covered by a ball ${\cal B}_d(Z,{ r })$)
 for the coverage area with   vol$\left([-1,1]^d \cap {\cal B}_d(\mathbb{Z}_n,r)\right)  $. Then, using accurate approximations for $C_{d,Z,{ r }}$, we derive   approximations for  vol$\left([-1,1]^d \cap {\cal B}_d(\mathbb{Z}_n,r)\right)  $.
In Theorem~\ref{Asymptotic_theorem} we derive asymptotic expressions for the $(1-\gamma)$-coverage radius for the design ${\mathbb{D}_{d,1/2}}$ and show that for any $\gamma>0$, the ratio of the $(1-\gamma)$-coverage radius to the $1$-coverage radius tends to $1/\sqrt{3}$ as $d \to \infty$. Numerical results of Section~\ref{sec:asymptotics} confirm that even for rather small $d$, the $0.999$-coverage radius is much smaller than the $1$-coverage radius providing the full coverage.

In Section~\ref{sec:numerical} we demonstrate that the approximations developed in Section~\ref{Main_results} are very accurate and make a comparative study of selected designs  used for quantization and covering.

In Appendices A--C, we provide proofs of the most technical results. In Appendix D, for completeness, we briefly derive an approximation for
$C_{d,Z,{ r }}$ with arbitrary $d$, $Z$ and $r$.

The two most important contributions  of this paper are: a) derivation of the closed-form expression for the quantization error for the design ${\mathbb{D}_{n,\delta}}$, and  b) derivation of accurate approximations for {the coverage area}  vol$\left([-1,1]^d \cap {\cal B}_d(\mathbb{Z}_n,r)\right) $
for the design ${\mathbb{D}_{n,\delta}}$.

\section{Quantization}

\label{sec:quant}

\subsection{Reformulation in terms of the Voronoi cells}

Consider any $n$-point design $\mathbb{Z}_n = \{Z_1, \ldots, Z_n\}$. The Voronoi cell $V(Z_i)$ for  $Z_i\in \mathbb{Z}_n$ is defined as
\bea
V(Z_i) = \{ x \in [-1,1]^d : \norm{Z_i-x} \leq \norm{Z_j-x} \text{ for } j\neq i \} \,.
\eea

\AZ{ The mean squared quantization error $\theta(\mathbb{Z}_n)$ introduced in \eqref{eq:errorQ} can be written in terms of the Voronoi cells as follows:
\be
\label{eq:Vor}
\theta(\mathbb{Z}_n)=\mathbb{E}_X\min_{i=1, \ldots, n} \|X-Z_i\|^2 = \frac{1}{{\rm vol ([-1,1]^d)}}\sum_{i=1}^{n} \int_{V(Z_i)}\norm{ X-Z_i }^2 dX \,,
\ee
where $ X = (x_1,  \ldots, x_d) $ and $dX = dx_1 dx_2\cdots dx_d$.}

This reformulation has significant benefit when the design $\mathbb{Z}_n$ has certain   structure.
In particular, if all  of the Voronoi cells $V(Z_i), i=1,\ldots,n, $ are congruent, then we can simplify \eqref{eq:Vor} to
\be
\label{eq:theta_simple}
\theta(\mathbb{Z}_n) = \frac{1}{{\rm vol}(V(Z_1))}\int_{V(Z_1)}\norm{ X-Z_1 }^2 dX \,.
\ee
In Section~\ref{quant_exact_results}, this formula will be the starting point for derivation of   the closed-form expression for $\theta(\mathbb{Z}_n)$  for the design $\mathbb{D}_{n, \delta}$.

\subsection{Re-normalization of the quantization error }
To compare  efficiency of $n$-point designs $\mathbb{Z}_n$  with different values of $n$, one must suitably normalise $\theta(\mathbb{Z}_n)$ with respect to $n$. Specialising a classical characteristic for quantization in space, as formulated in \cite[f-la (86), Ch.2]{Conway}, 
we obtain
\be
\label{eq:norm_G}
Q_d (\mathbb{Z}_n) = \frac{1}{d} \frac{\frac{1}{n}\sum_{i=1}^{n} \int_{V(Z_i)}\norm{ X-Z_i }^2 dX}{ \left[ \frac1n \sum_{i=1}^{n}
{\rm vol}(V(Z_i))\right]^{1+\frac2d} } \,.
\ee
Note that $Q_d (\mathbb{Z}_n)$ is re-normalised with respect to dimension $d$ too, not only with respect to $n$. Normalization $1/d$ with respect to $d$ is very natual in view of the definition of the Euclidean norm.

Using \eqref{eq:Vor}, for the cube $[-1,1]^d$, \eqref{eq:norm_G}  can be expressed as
\be\label{eq:norm_Qd}
Q_d (\mathbb{Z}_n) = \frac{n^{2/d}\theta(\mathbb{Z}_n)}{d\left[\sum_{i=1}^{n}
{\rm vol}(V(Z_i)) \right]^{2/d}} =  \frac{n^{2/d}\theta(\mathbb{Z}_n)}{d\cdot
{\rm vol}([-1,1]^d) ^{2/d}} = \frac{n^{2/d}}{4d} \theta(\mathbb{Z}_n) \, .
\ee

%
%

\subsection{Voronoi cells for $\mathbb{D}_{n,\delta}$}
\label{Voronoi}

\begin{proposition}\label{voronoi_prop}
Consider the design ${\mathbb{D}_{n,\delta}^{(0)}}$,
the collection of   $n=2^{d}$ points $(\pm \delta , \ldots, \pm \delta)$, $0< \delta< 1$.
The Voronoi cells  for this  design  are all congruent. The Voronoi cell for the point $ \boldsymbol{\delta}=(\delta,\delta,\ldots,\delta)$ is
the cube
\be
\label{v0}
 C_0 = \left\{ X\! = \!(x_1,\ldots,x_d)\! \in\! \mathbb{R}^d\! : \;0\leq x_i \leq 1 ,\; i=1,2,\ldots, d\right\}\, .
\ee

\end{proposition}
\label{pr1}
$\mathbf{Proof}$.
Consider the  Voronoi cells created by the
design ${\mathbb{D}_{n,\delta}^{(0)}}$ in the whole space $\mathbb{R}^d$.
 For the point $ \boldsymbol{\delta}=(\delta,\delta,\ldots,\delta)$,  the  Voronoi cell is clearly $\{X\! = \!(x_1,\ldots,x_d)\!: \, x_i \geq 0\}$. By intersecting this set with the cube $[-1,1]^d$ we obtain \eqref{v0}.
  \hfill $\Box$\\

\begin{theorem} \label{th:1}
The Voronoi cells of the design $\mathbb{D}_{n,\delta} =  \{Z_1, \ldots, Z_n\}$ are all congruent. The Voronoi cell for the point
$Z_1 = \boldsymbol{\delta}=(\delta,\delta,\ldots,\delta)\in \mathbb{R}^d$ is
\be
\label{eq:uj1}
V(Z_1) = C_0 \bigcup \left[ \bigcup_{j=1}^{d} U_j \right]
\ee
where $C_0$ is the cube \eqref{v0} and
\be
\label{eq:uj}
U_j = \left\{  X = (x_1,x_2,\ldots,x_d)\! \in\! \mathbb{R}^d\! : -1\leq x_j \leq 0 , |x_j|\leq x_k \leq 1\; \mbox{ \rm for all } k\neq j\right\} \, .
\ee
The volume of $V(Z_1)$ is ${{\rm vol}(V(Z_1))}=2$.

\end{theorem}

$\mathbf{Proof}$. The design $\mathbb{D}_{n,\delta}$ is symmetric with respect to all components implying that all $n=2^{d-1}$ Voronoi cells are congruent immediately yielding that  their volumes  equal 2. Consider $V(Z_1)$ with $Z_1 = \boldsymbol{\delta}$.

Since $\mathbb{D}_{n,\delta} \subset {\mathbb{D}_{n,\delta}^{(0)}}$,  where design ${\mathbb{D}_{n,\delta}^{(0)}}$ is introduced in Proposition \ref{voronoi_prop}, and $C_0 $ is the Voronoi set of $\boldsymbol{\delta}$ for design ${\mathbb{D}_{n,\delta}^{(0)}}$, $C_0 \subset V(\boldsymbol{\delta})$
for design $\mathbb{D}_{n,\delta}$ too.

Consider the $d$ cubes adjacent to $C_0$:
\be\label{adjacent_cubes}
C_j = \left\{  X = (x_1,x_2,\ldots,x_d)\! \in\! \mathbb{R}^d\! : -1\leq x_j \leq 0 ,\; 0\leq x_i \leq 1\; \mbox{ \rm for all } i\neq j\right\} \, ;\;\; j=1, \ldots, d.
\ee
A part of each cube $C_j$  belongs to $V(Z_1)$. This part is exactly the set $U_j$ defined by \eqref{eq:uj}. This can be seen as follows.
A part of $C_j$  also belongs to the Voronoi set of the point $X_{jk}= \boldsymbol{\delta}-2\delta e_j-2\delta e_k$, where $e_l=(0,\ldots,0,1,0,\ldots,0)$ with 1 placed at $l$-th place; all components of $X_{jk}$ are $\delta$ except $j$-th and $k$-th components which are $-\delta$. We have to have $|x_j| \leq x_k$, for a point $X \in C_j$ to be closer to $Z_1$ than to $X_{jk}$. Joining all constraints for
$ X = (x_1,x_2,\ldots,x_d) \in C_j$ ($k=1, \ldots, d$, $k \neq j$) we obtain \eqref{eq:uj} and hence \eqref{eq:uj1}.

  \hfill $\Box$\\

%
%

\subsection{Explicit formulae for the quantization error}
\label{quant_exact_results}

\begin{theorem}\label{quant_theorem}
For the design $\mathbb{D}_{n,\delta}$ with $0 \leq \delta\leq 1$, we obtain:
\be
\theta(\mathbb{D}_{n,\delta}) &=& d\left(\delta^2-\delta + \frac13 \right) + \frac{2\delta}{d+1}\, , \label{th}\\
Q_d (\mathbb{D}_{n,\delta})&=& 2^{-2/d}\left(  \delta^2-\delta + \frac13 + \frac{2\delta}{d(d+1)} \right) \label{thQ} \, .
\ee
\end{theorem}
$\mathbf{Proof}$.
To compute $\theta(\mathbb{D}_{n,\delta})$, we use \eqref{eq:theta_simple}, where, in view of Theorem~\ref{th:1}, ${{\rm vol}(V(Z_1))}=2$.
Using the expression \eqref{eq:uj1} for $V(Z_1)$ with  $Z_1 = \boldsymbol{\delta}$, we obtain
\be
\theta(\mathbb{Z}_n) = \frac{1}{2}\int_{V(Z_1)}\norm{ X-Z_1 }^2 dX = \frac{1}{2}\left[\int_{C_0}\norm{ X-Z_1 }^2 dX + d \int_{U_1}\norm{ X-Z_1 }^2 dX  \right]\, . \label{twoterms}
\ee
Consider the two terms in \eqref{twoterms} separately. The first term is easy:
\be \label{useful_integral}\int_{C_0}\norm{ X-Z_1 }^2 dX = \int_{C_0} \sum_{i=1}^{d} (x_i-\delta)^2 dx_1\ldots dx_d = d \int_{0}^{1} (x-\delta)^2 dx= d\left(\delta^2-\delta+\frac13\right)\, .
\ee
For the second term we have:
\be
\!\!\!\!\!\!\int_{U_1}\norm{ X-Z_1 }^2 dX & =& \int_{-1}^{0} \left[ \int_{|x_1|}^{1} \ldots \int_{|x_1|}^{1} \sum_{i=1}^{d} (x_i-\delta)^2 dx_2 \ldots dx_d  \right] dx_1 \nonumber\\
& =&  \int_{-1}^{0} (x_1-\delta)^2 (1+x_1)^{d-1} dx_1 + (d-1) \int_{-1}^{0} (1+x_1)^{d-2} \int_{|x_1|}^{1} (x_2-\delta)^2 dx_2 dx_1 \nonumber\\
& =& \delta^2-\delta+\frac13 + \frac{4 \delta}{d(d+1)}\, . \label{int_in_seq}
\ee
Inserting the obtained expressions into \eqref{twoterms} we obtain \eqref{th}. The expression \eqref{thQ} is a consequence of \eqref{eq:norm_Qd}, \eqref{th} and $n=2^{d-1}$.
 \hfill $\Box$\\

A simple consequence of Theorem~\ref{quant_theorem} is the following corollary.
\begin{corollary} \label{cor:1}
The optimal value of $\delta$ minimising $\theta(\mathbb{D}_{n,\delta})$ and $Q_d (\mathbb{D}_{n,\delta})$ is
\be
\label{eq:opt_d}
\delta^\ast = \frac12 - \frac{1}{d(d+1)} \,;
\ee
for this value,
\be
\label{eq:opt_d1}
Q_d (\mathbb{D}_{n,\delta^\ast})= \min_{\delta} Q_d (\mathbb{D}_{n,\delta})=2^{-2/d}\left[  \frac1{12}+{\frac {{d}^{2}+d-1}{ \left( d+1 \right) ^{2}{d}^{2}}} \right]\, .
\ee
\end{corollary}
Let us make several remarks.

\begin{enumerate}

\item

\AZ{The value $\delta^\ast$ can be alternatively characterised by the well-known optimality condition of a general design saying that each design point of an optimal quantizer must be a centroid of the related Voronoi cell; see e.g. \cite{saka2007latinized}. Specifically, each design points $Z_i \in {\mathbb{D}_{n,\delta}}$ is the  centroid of
$V(Z_i)$ if and only if $\delta=\delta^\ast$.}

\item
From \eqref{thQ},  for the design $\mathbb{D}_{n,1/2}$ we get
\be
\label{eq:opt_d0}
Q_d (\mathbb{D}_{n,1/2})= 2^{-2/d}\left[  \frac1{12}+\frac {1}{ \left( d+1 \right) d} \right]\,;
\ee
this value is always slightly larger than \eqref{eq:opt_d1}.

\item
For the one-point design $\mathbb{D}^{(0)}= \{0\}$ with  the single point 0 and the design $\mathbb{D}_{n}^{(0)}$ with $n=2^d$ points $(\pm \frac12, \ldots, \pm \frac12)$ we have $Q_d (\mathbb{D}^{(0)})=Q_d (\mathbb{D}_{n}^{(0)})=1/12$, which coincides with the value of $Q_d $ in the case of space quantization by the integer-point lattice $Z^d$, see \cite[Ch. 2 and 21]{Conway}.

  \item
The quantization error \eqref{eq:opt_d0} for  the design ${\mathbb{D}_{n,1/2}}$ have almost exactly the same form as the
 quantization error for the `checkerboard lattice' $D_d$ in $\mathbb{R}^d$; the difference is in the factor $1/2$ in the last term in  \eqref{eq:opt_d0}, see \cite[f-la (27), Ch.21]{Conway}. Naturally, the quantization error $Q_d$ for $D_d$ in $\mathbb{R}^d$ is slightly smaller than $Q_d$ for $\mathbb{D}_{n,1/2}$ in $[-1,1]^d$.

  \item The optimal value of $\delta$ in \eqref{eq:opt_d} is smaller than $1/2$. This is caused by a non-symmetrical shape of the Voronoi cells $V(Z_j)$ for designs $ \mathbb{D}_{n,\delta}$, which is clearly visible in \eqref{eq:uj1}.

\item The minimal value of $Q_d (\mathbb{D}_{n,\delta^\ast})$ with respect to $d$ is attained at $d=15$.

  \item Formulas \eqref{eq:opt_d} and \eqref{eq:opt_d1} are in agreement  with  numerical results presented in Table~4 of \cite{us} and Table~5 of~\cite{second_paper}.

\end{enumerate}

\begin{figure}[h]
\centering
\begin{minipage}{.5\textwidth}
  \centering
  \includegraphics[width=1\textwidth]{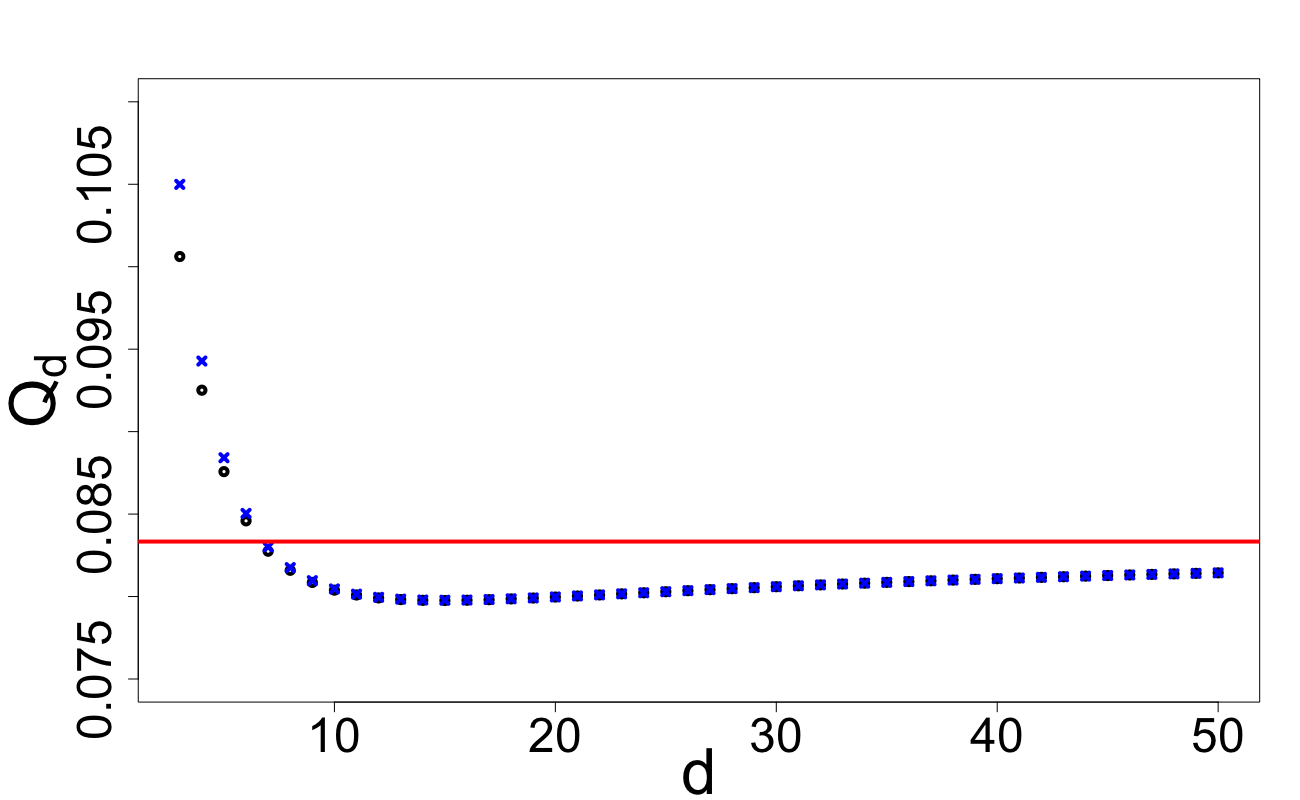}
\caption{$Q_d (\mathbb{D}_{n,\delta^\ast}) $  and $Q_d (\mathbb{D}_{n,1/2}) $ as functions \\ of $d$ and $Q_d (\mathbb{D}_{n}^{(0)})=1/12$; $d=3, \ldots,50$. }
\label{Quant_example1}
\end{minipage}%
\begin{minipage}{.5\textwidth}
  \centering
 \includegraphics[width=1\textwidth]{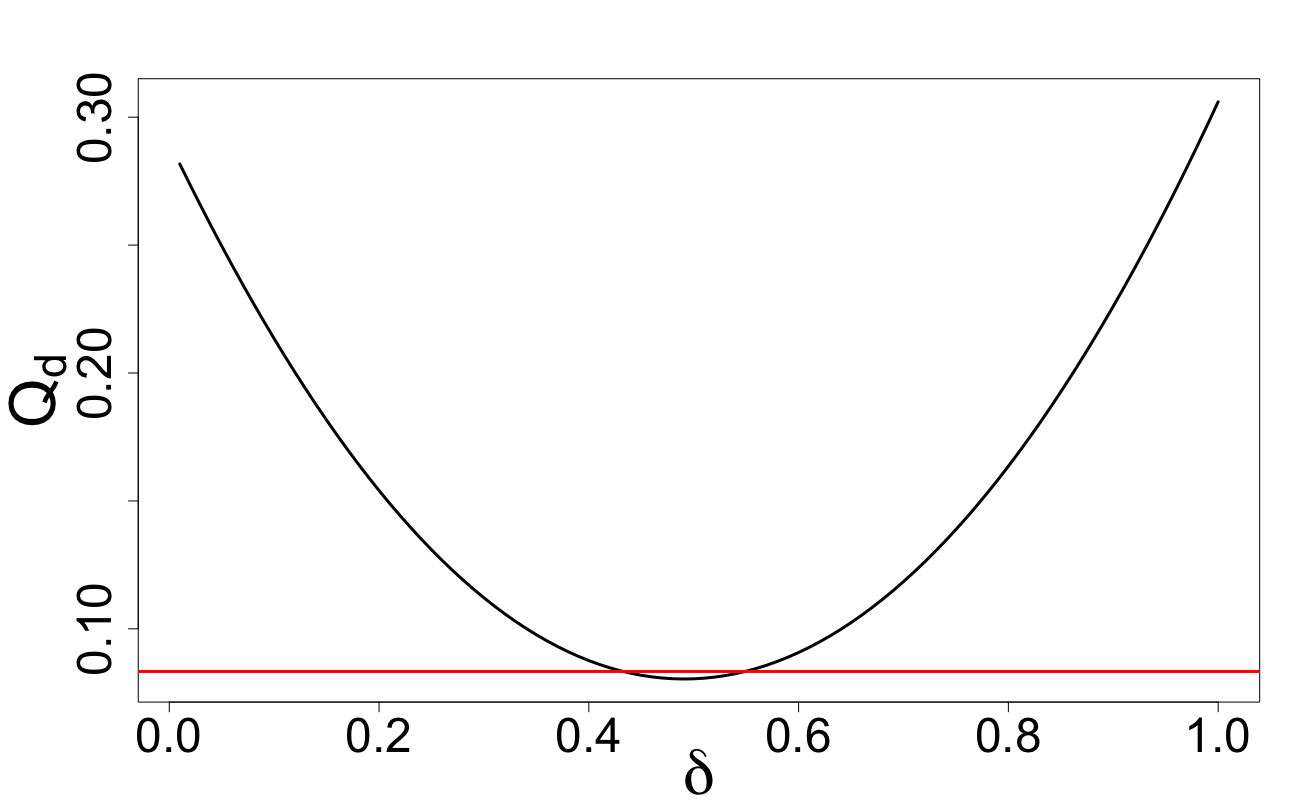}
\caption{$Q_d (\mathbb{D}_{n,\delta}) $ as a function of $\delta$ and \\$Q_d (\mathbb{D}_{n}^{(0)})=1/12$; $d=10$.}
\label{Quant_example2}
\end{minipage}
\end{figure}

Let us now briefly illustrate the results above. In Figure~\ref{Quant_example1}, the black circles depict the quantity $Q_d (\mathbb{D}_{n,\delta^\ast})$ as a function of $d$. The quantity $Q_d (\mathbb{D}_{n}^{(0)})=1/12$ is shown with the solid red line. We conclude that from dimension seven onwards, the design $\mathbb{D}_{n,\delta^\ast}$ provides better  quantization per points than the design $\mathbb{D}_{n}^{(0)}$. Moreover for $d> 15$, the quantity $Q_d (\mathbb{D}_{n,\delta^\ast})$ slowly increases and converges to $1/12$. Typical behaviour of $Q_d (\mathbb{D}_{n,\delta})$ as a function of $\delta$ is shown in Figure~\ref{Quant_example2}. This figure demonstrates the significance of choosing $\delta$ optimally.

\section{Closed-form expressions for the coverage area with  $\mathbb{D}_{n,\delta}$ and approximations}\label{Main_results}
In this section, we will derive explicit expressions for  the coverage area of the cube $[-1,1]^d$ by the union of the balls
${\cal B}_d({\mathbb{D}_{n,\delta}},r)$ associated with the   design ${\mathbb{D}_{n,\delta}}$ introduced in Section~\ref{main}.
That is, we will derive expressions for the quantity $C_d(\mathbb{D}_{n,\delta},r)$ for all values of $r$. Then, in Section~\ref{sec:app1}, we shall obtain approximations for $C_d(\mathbb{D}_{n,\delta},r)$. The accuracy of the approximations will be assessed in
Section~\ref{sec:accur1}.

\subsection{Reduction to Voronoi cells}
For an $n$-point design {$\mathbb{Z}_n = \{Z_1, \ldots, Z_n\}$}, denote the proportion of the Voronoi cell around $Z_i$ covered by the ball ${\cal B}_d(Z_i,{ r })$ as
\bea
V_{d,Z_i,{ r }}:=  {\rm vol}(V(Z_i) \cap  {\cal B}_d(Z_i,{ r }))/{\rm vol}(V(Z_i))\,  .
\eea
The following  lemma is straightforward.

\begin{lemma} \label{lem:equal}
Consider a design  $\mathbb{Z}_n=\{Z_1, \ldots, Z_n\}$ such that all Voronoi cells $V(Z_i)$ are congruent. Then for any
 $Z_i \in \mathbb{Z}_n$, $ C_d(\mathbb{Z}_n,r)=V_{d,Z_i,{ r }}$.
\end{lemma}

In view of Theorem~\ref{th:1}, for design $\mathbb{D}_{n,\delta}$  all Voronoi cells $V(Z_i)$ are congruent and ${\rm vol}(V(Z_i))=2$; recall that $n=2^{d-1}$.
By then applying Lemma~\ref{lem:equal} and without loss of generality we have choosen $Z_1 = \boldsymbol{\delta}= (\delta,\delta,\ldots,\delta)\in \mathbb{R}^d$, we have for any $r>0$
\be
\label{voroni_ratio}
V_{d,\boldsymbol{\delta} ,{ r }}= \frac12 {\rm vol}(V(\boldsymbol{\delta}) \cap  {\cal B}_d(\boldsymbol{\delta},{ r })) =C_d(\mathbb{D}_{n,\delta},r)  \, .
\ee

 In order to formulate explicit expressions for $V_{d,\boldsymbol{\delta} ,{ r }}$, we need  the  important quantity, proportion of intersection of $[-1,1]^d$ with one ball.
Take the cube $[-1,1]^d$ and a ball
$
{\cal B}_d(Z,{ r })= \{ Y \in \mathbb{R}^d: \| Y-Z \| \leq { r } \}
$  centered at a point  $Z=(z_1, \ldots, z_d) \in \mathbb{R}^d$; this point $Z$ could be outside  $ [-1,1]^d$.
The fraction of the cube $[-1,1]^d$ covered by the ball ${\cal B}_d(Z,{ r })$ is denoted by
\bea
C_{d,Z,{ r }}={{\rm vol}([-1,1]^d \cap  {\cal B}_d(Z,{ r }))}/2^d\, .
\eea

%
%

\subsection{Expressing $C_d(\mathbb{D}_{n,\delta},r)$ through $C_{d,Z,{ r }}$ }

\begin{theorem}\label{Main_theorem}
Depending on the values of $r$ and $\delta$, the quantity $C_d(\mathbb{D}_{n,\delta},r)$ can be expressed through $C_{d,Z,{ r }}$ for suitable $Z$ as follows.
\begin{itemize}
\item For $r\leq \delta$:
\be \label{Exact_1}
C_d(\mathbb{D}_{n,\delta},r)=\frac12 C_{d,\boldsymbol{2\delta - 1},{2r }} \,.
\ee
\item For $\delta \leq r \leq 1+\delta$:
\be \label{Exact_2}
C_d(\mathbb{D}_{n,\delta},r)
&=& \frac12 \left[ C_{d,\boldsymbol{2\delta - 1},{2r }} + d \int_{0}^{r-\delta} C_{d-1,\boldsymbol{\frac{2\delta-1-{t} }{1-{t}}},{ \frac{2\sqrt{r^2-(t+\delta)^2}}{1-t} }} (1-t)^{d-1}\, dt  \right] \,.
\ee
\item For $r \geq 1+\delta$:
\be\label{Exact_3}
C_d(\mathbb{D}_{n,\delta},r) = \frac12 \left[ C_{d,\boldsymbol{2\delta - 1},{2r }} + d \int_{0}^{1} C_{d-1,\boldsymbol{\frac{2\delta-1- {t} }{1-{t}}},{ \frac{2\sqrt{r^2-(t+\delta)^2}}{1-t} }} (1-t)^{d-1}\, dt  \right] \,.
\ee

\end{itemize}

\end{theorem}

The proof of Theorem~\ref{Main_theorem} is given in Appendix A.

\subsection{Approximation for $C_d(\mathbb{D}_{n,\delta},r)$}
\label{sec:app1}

Accurate approximations for $C_{d,Z,{ r }}$ for arbitrary $d, Z$ and $r$ were developed in \cite{us}. By using the general expansion in the central limit theorem for sums of independent non-identical r.v., the following approximation was developed:
\be
\label{eq:inters2f_corrected}
C_{d,Z,{ r }} \cong \Phi(t) + \frac{  \|Z\|^2+d/63}{5\sqrt{3} (\|Z\|^2+d/15)^{3/2} }(1-t^2)\varphi(t) \, ,
\ee
where
\bea
t =
 \frac{\sqrt{3}(r^2- \|Z\|^2 -d/3)}{2\sqrt{ \|Z\|^2 +{d}/{15}} }\,.\\\nonumber
\eea
A short derivation of this approximation is included in Appendix D. 
Using \eqref{eq:inters2f_corrected}, we formulate the following approximation for $C_d(\mathbb{D}_{n,\delta},r)$. \\

{\bf Approximations for $C_d(\mathbb{D}_{n,\delta},r)$.} {\it   Approximate the values  $C_{\cdot,\cdot,\cdot}$ in formulas \eqref{Exact_1},\eqref{Exact_2},\eqref{Exact_3} with corresponding approximations  \eqref{eq:inters2f_corrected}.
}

\subsection{Simple bounds for $C_d(\mathbb{D}_{n,\delta},r)$}
\begin{lemma}\label{Upper_and_lower_bounds}
 For any $r \geq 0$, $ 0<\delta <1 $ and $\boldsymbol{\delta}= (\delta,\delta,\ldots,\delta) \in \mathbb{R}^d$, the quantity $C_d(\mathbb{D}_{n,\delta},r)$ can be bounded as follows:
\be
\label{eq:bounds}
\frac{1}{2} [  C_{d,\boldsymbol{2\delta-1},{2r }} +C_{d,A,{2r }}  ]  \leq C_d(\mathbb{D}_{n,\delta},r)  \leq C_{d,\boldsymbol{2\delta-1},{2r }} \,.
\ee
where $A= \left(2\delta+1,2\delta-1,\ldots,2\delta-1 \right) \in \mathbb{R}^d $.
\end{lemma}
The proof of Lemma~\ref{Upper_and_lower_bounds} is given in Appendix B.

In Figures~\ref{Lower1} and \ref{Lower2}, using the approximation given in~\eqref{eq:inters2f_corrected} we study the tightness of the bounds given in~\eqref{eq:bounds}. In these figures, the dashed red line, dashed blue line and solid black line depict the upper bound, the lower bound and the approximation for $C_d(\mathbb{D}_{n,\delta},r)$ respectively. We see that the upper bound is very sharp across $r$ and $d$; this behaviour is not seen with the lower bound.
\begin{figure}[h]
\centering
\begin{minipage}{.5\textwidth}
  \centering
  \includegraphics[width=1\textwidth]{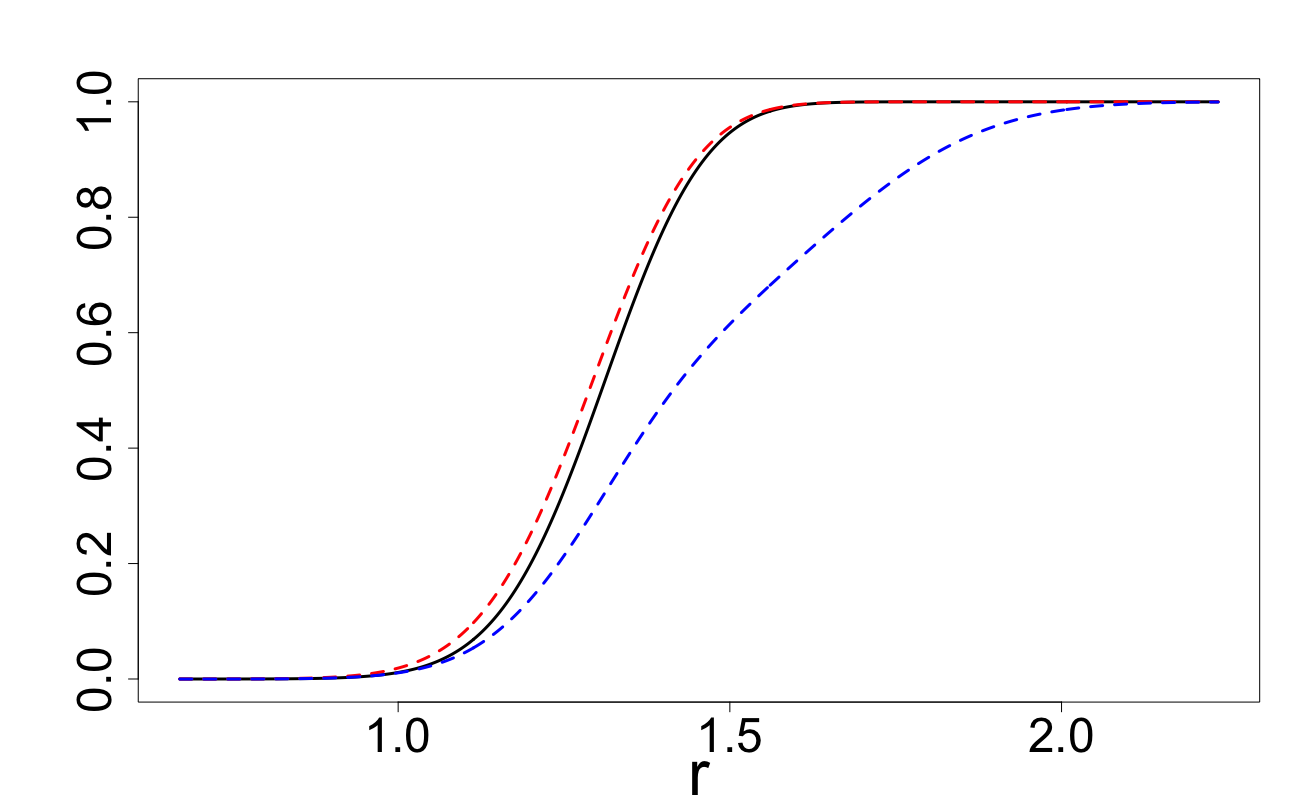}
\caption{$C_d(\mathbb{D}_{n,\delta},r)$ with upper and lower bounds:\\ $d=20$.}
\label{Lower1}
\end{minipage}%
\begin{minipage}{.5\textwidth}
  \centering
  \includegraphics[width=1\textwidth]{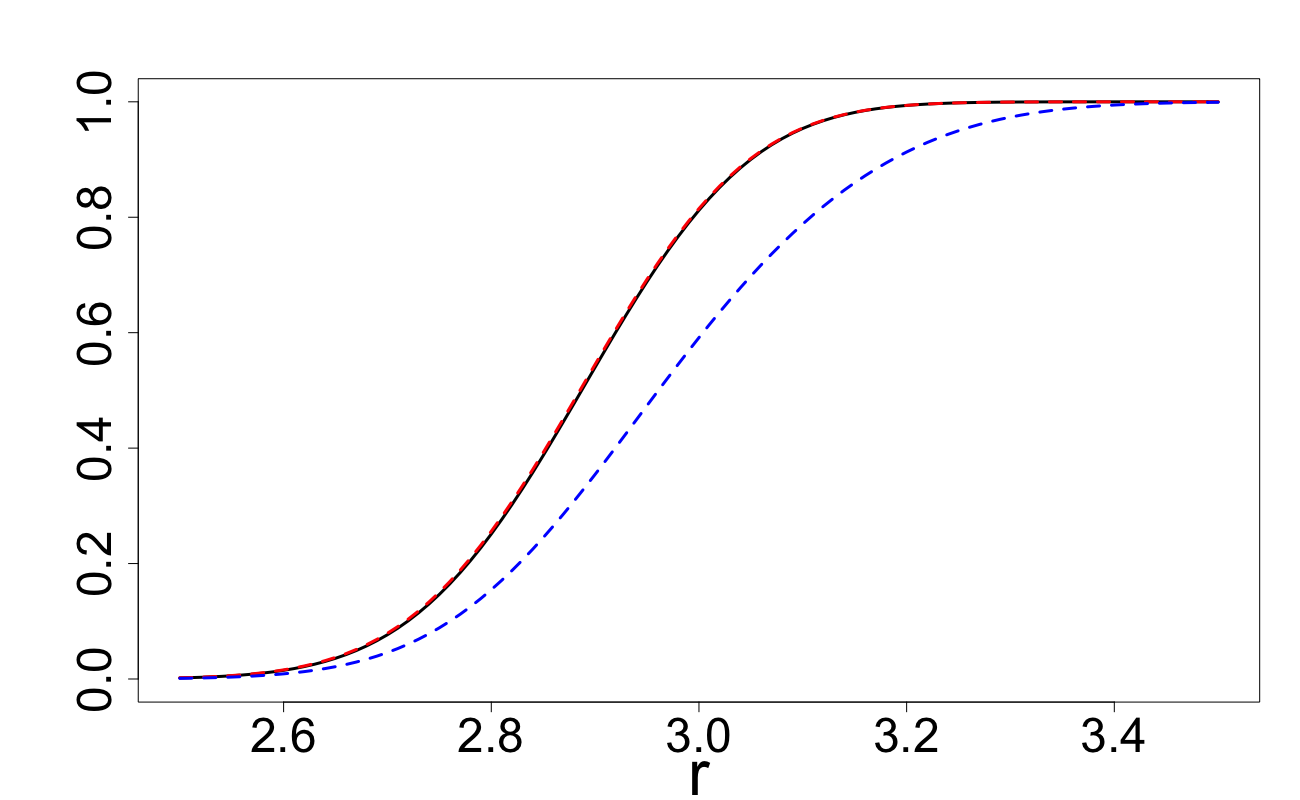}
\caption{$C_d(\mathbb{D}_{n,\delta},r)$ with upper and lower bounds: \\$d=100$.}
\label{Lower2}
\end{minipage}
\end{figure}

\subsection{`Do not try to cover the vertices'}
\label{sec:asymptotics}

In this section, we theoretically support the recommendation `do not try to cover the vertices' which was first stated in \cite{us} and supported in \cite{second_paper} on the basis of numerical evidence. In other words, we will show on the example of the design $\mathbb{D}_{n,1/2}$  that in large dimensions the attempt to cover the whole cube rather than 0.999 of it leads to a dramatic increase of the required radius of the balls.

\begin{theorem}\label{Asymptotic_theorem}
Let $\gamma$ be fixed, $0\leq \gamma\leq 1$. Consider $(1-\gamma)$-coverings  of $[-1,1]^d$  generated by the designs ${\mathbb{D}_{n,\delta}}$ and the associated normalized radii
$ R_{1\!-\!\gamma} ({\mathbb{D}_{n,\delta}})$, see \eqref{eq:R_g}.
For any $0<\gamma<1$ and  $0\leq \delta \leq 1$, the limit of $ R_{1\!-\!\gamma} ({\mathbb{D}_{n,\delta}})$, as $d\rightarrow \infty$, exists and
achieves minimal value  for  $\delta=1/2$. Moreover,
 $ R_{1\!-\!\gamma} ({\mathbb{D}_{n,1/2}})/R_{1} ({\mathbb{D}_{n,1/2}}) \rightarrow 1/\sqrt{3}$ as $d\rightarrow \infty,$ for  any $0<\gamma<1$.
\end{theorem}

Proof is given in Appendix C.

\AZ{ In Figures~\ref{Normalised_1}-\ref{Normalised_2} using a solid red line we depict the approximation of $C_d({\mathbb{D}_{n,1/2}},r)$ as a function of $R ={n^{1/d}} r/({ 2 \sqrt{d}} )$, see   \eqref{eq:Rr}. The vertical green line illustrates the value of $R_{0.999}$ and the vertical blue line depicts $R_{1} = {n^{1/d}} \sqrt{d+8}/({ 4 \sqrt{d}} )$. These figures illustrate that as $d$ increases, for all $\gamma$ we have $R_{1-\gamma}/R_{1}$ slowly tending to $1/\sqrt{3}$. From the proof of Theorem~\ref{Asymptotic_theorem}, it transpires that $C_d({\mathbb{D}_{n,\delta}},r)$ as a
function of $R$ converges to the jump function with the jump at  $1/(2\sqrt{3})$.}

%

\begin{figure}[H]
\centering
\begin{minipage}{.5\textwidth}
  \centering
  \includegraphics[width=1\textwidth]{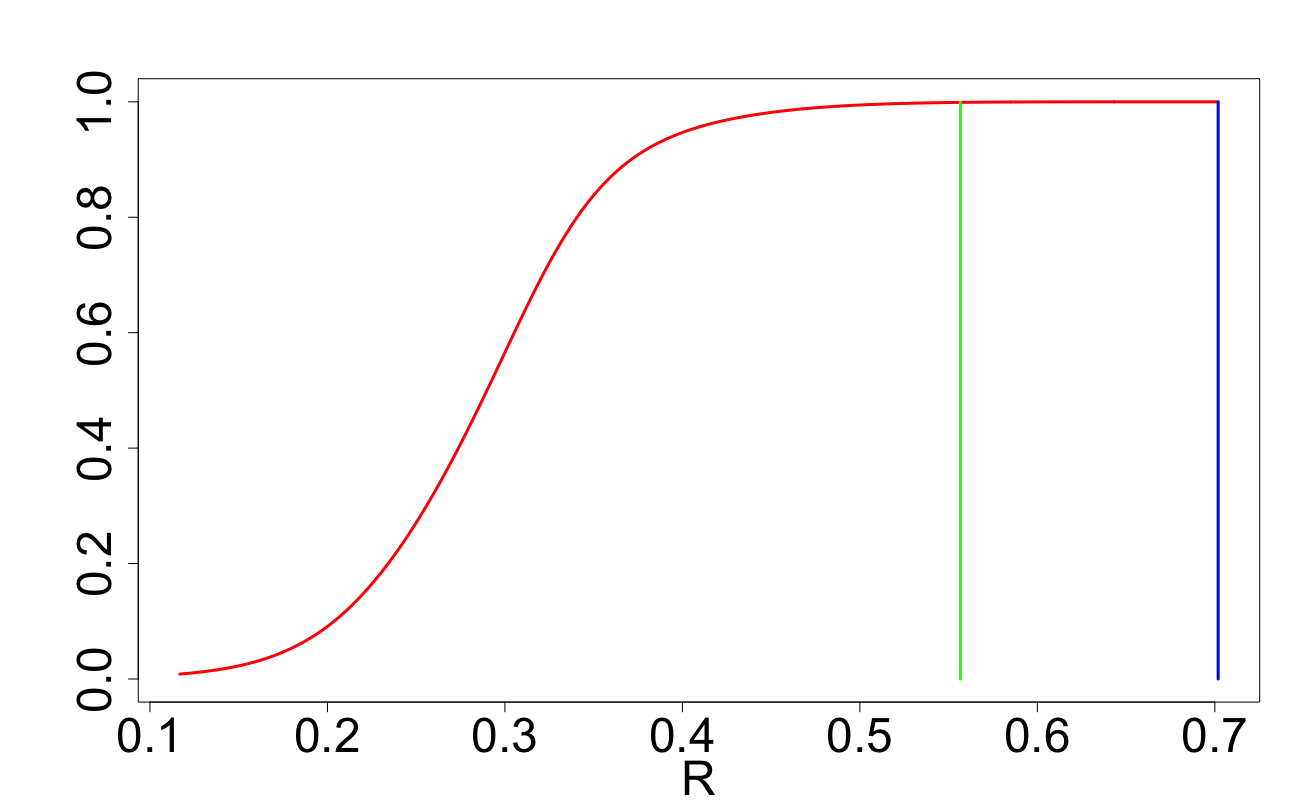}
\caption{$C_d({\mathbb{D}_{n,1/2}},r_{})$ with $R_{0.999}$ and $R_{1}$: $d=5$.}
\label{Normalised_1}
\end{minipage}%
\begin{minipage}{.5\textwidth}
  \centering
  \includegraphics[width=1\textwidth]{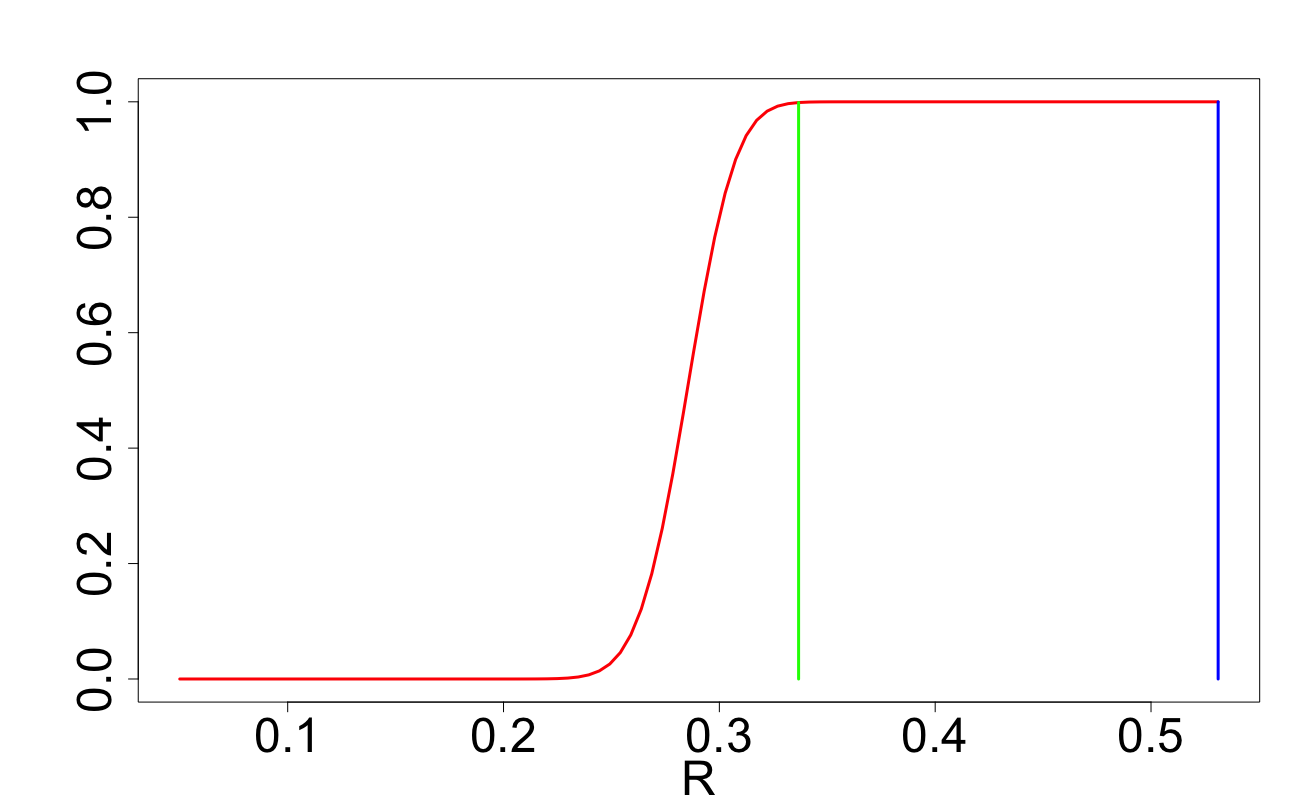}
\caption{$C_d({\mathbb{D}_{n,1/2}},r_{})$ with $R_{0.999}$ and $R_{1}$: $d=50$.}
\label{Normalised_2}
\end{minipage}
\end{figure}

%

\section{Numerical studies}
\label{sec:numerical}

For comparative purposes, we introduce another  design which is one of the most popular designs (both, for quantization and covering) considered in applications.\\

{\bf Design $\mathbb{S}_{n}$:} {\it
$Z_1, \ldots, Z_n$ are taken from a low-discrepancy Sobol's sequence on  the cube  $[-1,1]^d$.\\
}

For constructing the design $\mathbb{S}_{n}$, we use the R-implementation provided in the well-known `SobolSequence' package  \cite{Sobol}.
For $\mathbb{S}_{n}$, we have set $n=1024$ and $F2=10$ (an input parameter for the Sobol sequence function). \AZ{Sobol sequences $\mathbb{S}_{n}$
attain their best  space-filling properties when  $n$ is a power of $2$; that is, when $n=2^\ell$ for some integer $\ell$. We have chosen $\ell=10$.  As we study renormalised characteristics $Q_d (\cdot)$ and $R_{1\!-\!\gamma} (\cdot)$ of designs, exact value of $\ell$  for $\mathbb{S}_{n}$ with $n=2^\ell$ is almost irrelevant: in particular, numerically computed values $Q_d (\mathbb{S}_{2^\ell})$ and $R_{1\!-\!\gamma} (\mathbb{S}_{2^\ell})$ for $\ell =8,9,11,12$ are almost indistinguishable from the corresponding values for $\ell=10$ provided below in Tables~\ref{Quant_table} and ~\ref{Covering_table}.  By varying values of $\ell$, we are not improving space-filling properties of  $\mathbb{S}_{2^\ell}$. In fact, increase of $\ell$ generally leads to a slight deterioration of normalised space-filling characteristics (including $Q_d (\cdot)$ and $R_{1\!-\!\gamma} (\cdot)$) of Sobol sequences. }

\subsection{Quantization and weak covering comparisons}
\label{sec:compar}

In Table~\ref{Quant_table}, we compare the normalised mean squared quantization error $Q_d (\mathbb{Z}_n)$ defined in \eqref{eq:norm_Qd} across three designs: $\mathbb{D}_{n,\delta^\ast}$ with $\delta^\ast$ given in \eqref{eq:opt_d}, $\mathbb{D}_{n}^{(0)}$ and $\mathbb{S}_{n}$.

\begin{table}[h]
\centering
\begin{tabular}{ |p{2.7cm}||p{2cm}|p{2cm}|p{2cm}|p{2cm}|p{2cm}|  }

 \hline
  & $d=5$ & $d=7$&$d=10$& $d=15$ & $d=20$\\
 \hline
$Q_d (\mathbb{D}_{n,\delta^\ast})$&  0.0876   & 0.0827  & 0.0804  &  0.0798 &  0.0800\\
$Q_d (\mathbb{D}_{n}^{(0)})$ &   0.0833   & 0.0833   & 0.0833 & 0.0833 & 0.0833  \\
$Q_d (\mathbb{S}_{n})$  & 0.0988    & 0.1003 &0.1022    & 0.1060   &  0.1086\\
 \hline
\end{tabular}
\caption{Normalised mean squared quantization error $Q_d$ for three designs and different $d$. }
\label{Quant_table}
\end{table}

In Table~\ref{Covering_table}, we compare the normalised statistic $R_{1\!-\!\gamma}$ introduced in \eqref{eq:thick4}, where we have fixed $\gamma = 0.01$. For designs $\mathbb{D}_{n,\delta}$ (with the optimal value of $\delta$), $ \mathbb{D}_{n,{\tiny 1/2}}$ and $\mathbb{D}_{n}^{(0)}$ we have also included $R_{1}$, the smallest normalised radius that ensures the full coverage.

\begin{table}[h]
\begin{center}
\begin{tabular}{ |p{2.7cm}||p{2cm}|p{2cm}|p{2cm}|p{2cm}|p{2cm}|  }
 \hline
  & $d=5$ & $d=7$&$d=10$& $d=15$ & $d=20$\\
 \hline
$R_{1\!-\!\gamma}(\mathbb{D}_{n,\delta})$&   0.4750 (0.54)    & 0.3992 (0.53)  &  0.3635 (0.52)  & 0.3483  (0.51)  &  0.3417  (0.50)   \\
$R_{1\!-\!\gamma}(\mathbb{D}_{n,1/2})$ &  0.4765    & 0.4039   &   0.3649  &  0.3484  &  0.3417   \\
$R_{1\!-\!\gamma}(\mathbb{D}_{n}^{(0)})$ &  0.4092   & 0.3923  & 0.3766  & 0.3612 & 0.3522   \\
$R_{1\!-\!\gamma} (\mathbb{S}_{n})$  &  0.4714  &  0.4528  &  0.4256   &  0.4074  & 0.3967  \\
 \hline
\end{tabular}
\end{center}
\begin{center}
\begin{tabular}{ |p{2.7cm}||p{2cm}|p{2cm}|p{2cm}|p{2cm}|p{2cm}|  }
 \hline
$R_{1}(\mathbb{D}_{n,\delta})$&  0.6984  (0.54)   & 0.6555 (0.53) & 0.6178 (0.52)   &   0.5856 (0.51) &  0.5714  (0.50)   \\
$R_{1}(\mathbb{D}_{n,1/2})$ & 0.7019    &  0.6629  &  0.6259  &  0.5912   &   0.5714   \\
$R_{1}(\mathbb{D}_{n}^{(0)})$ & 0.5000   &  0.5000 & 0.5000  & 0.5000 & 0.5000    \\
 \hline
\end{tabular}
\end{center}
\caption{Normalised statistic $R_{1\!-\!\gamma}$ across $d$ with $\gamma=0.01$ (value in brackets corresponds to optimal $\delta$)}
\label{Covering_table}
\end{table}

Let us make some remarks concerning Tables~\ref{Quant_table} and \ref{Covering_table}.
\begin{itemize}

\item In conjunction with Figure~\ref{Quant_example1}, Table~\ref{Quant_table} shows that for $d\ge7$, the quantization for design $\mathbb{D}_{n,\delta^\ast}$ is superior over all other designs considered.
\item For the weak coverage statistic $R_{1\!-\!\gamma}$, the superiority of $\mathbb{D}_{n,\delta}$ with optimal $\delta$
over all other designs considered is seen for $d\ge 10$.
\item For the designs $\mathbb{D}_{n,\delta}$, the optimal value of $\delta$ minimizing  $R_{1\!-\!\gamma}$ depends on $\gamma$.
\item From remark 6 of  Section~\ref{quant_exact_results}, the minimal value of $Q_d (\mathbb{D}_{n,\delta^\ast})$ with respect to $d$ is attained at $d=15$.  For $d>15$, the quantity $Q_d (\mathbb{D}_{n,\delta^\ast})$ increases with $d$, slowly converging to $Q_d (\mathbb{D}_{n}^{(0)}) = 1/12$. This non-monotonic behaviour can be seen in  Table~\ref{Quant_table}.
\item Unlike the case of $Q_d (\mathbb{D}_{n,\delta^\ast})$, such non-monotonic behaviour is not seen for the quantity $R_{1\!-\!\gamma}$ and $R_{1\!-\!\gamma}(\mathbb{D}_{n,\delta})$ monotonically decreases  as $d$ increases. Also,  Theorem~\ref{Asymptotic_theorem} implies that
    for any $\gamma \in (0,1)$, $R_{1\!-\!\gamma}(\mathbb{D}_{n,\delta})\rightarrow {1}/({2\sqrt{3}})\cong 0.289$ as $d\rightarrow \infty$.

\end{itemize}

\subsection{Accuracy of covering approximation and dependence on $\delta$}

\label{sec:accur1}
\begin{figure}[h]
\centering
\begin{minipage}{.5\textwidth}
  \centering
 \includegraphics[width=1\textwidth]{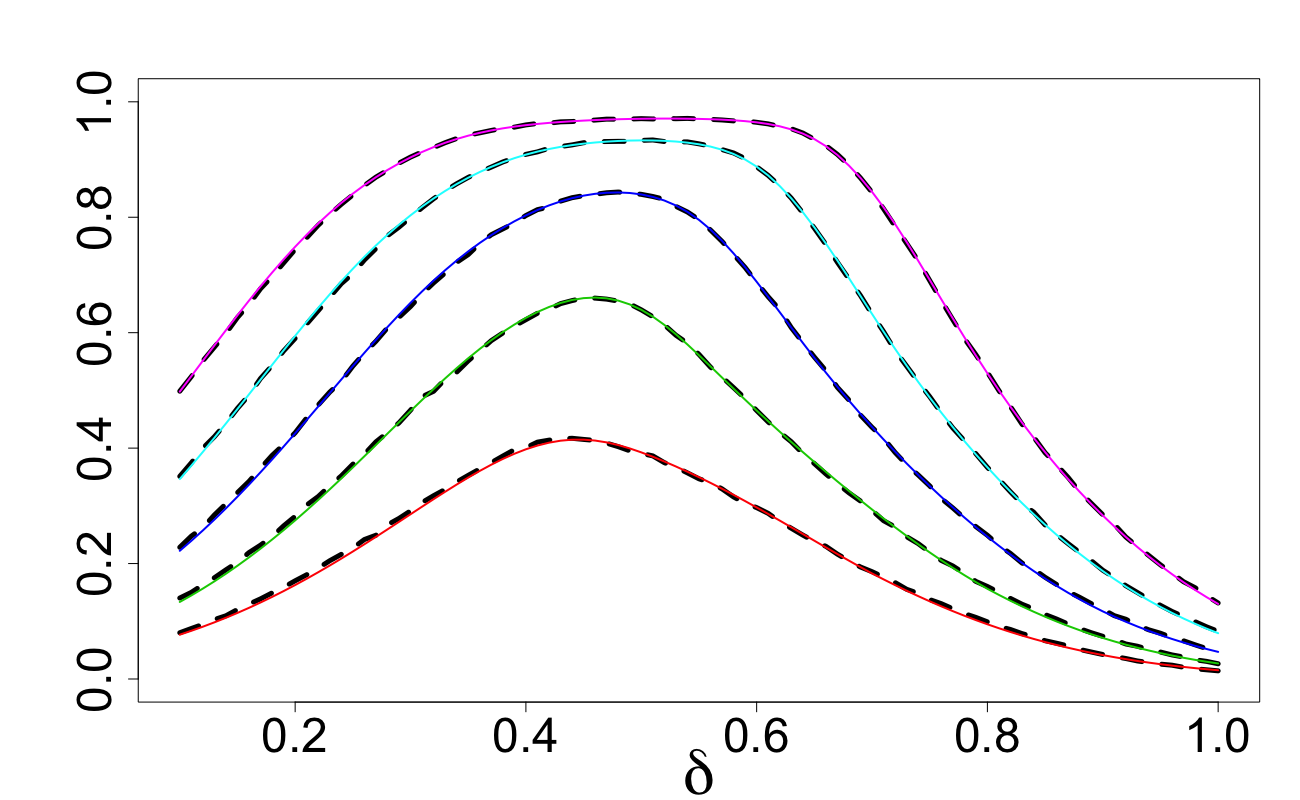}
\caption{$C_d(\mathbb{D}_{n,\delta},r)$ and its approximation: $d=5$, $\\r$ from $0.7$ to $1.1$ increasing by $0.1$}
\label{Covering_approx1}
\end{minipage}%
\begin{minipage}{.5\textwidth}
  \centering
 \includegraphics[width=1\textwidth]{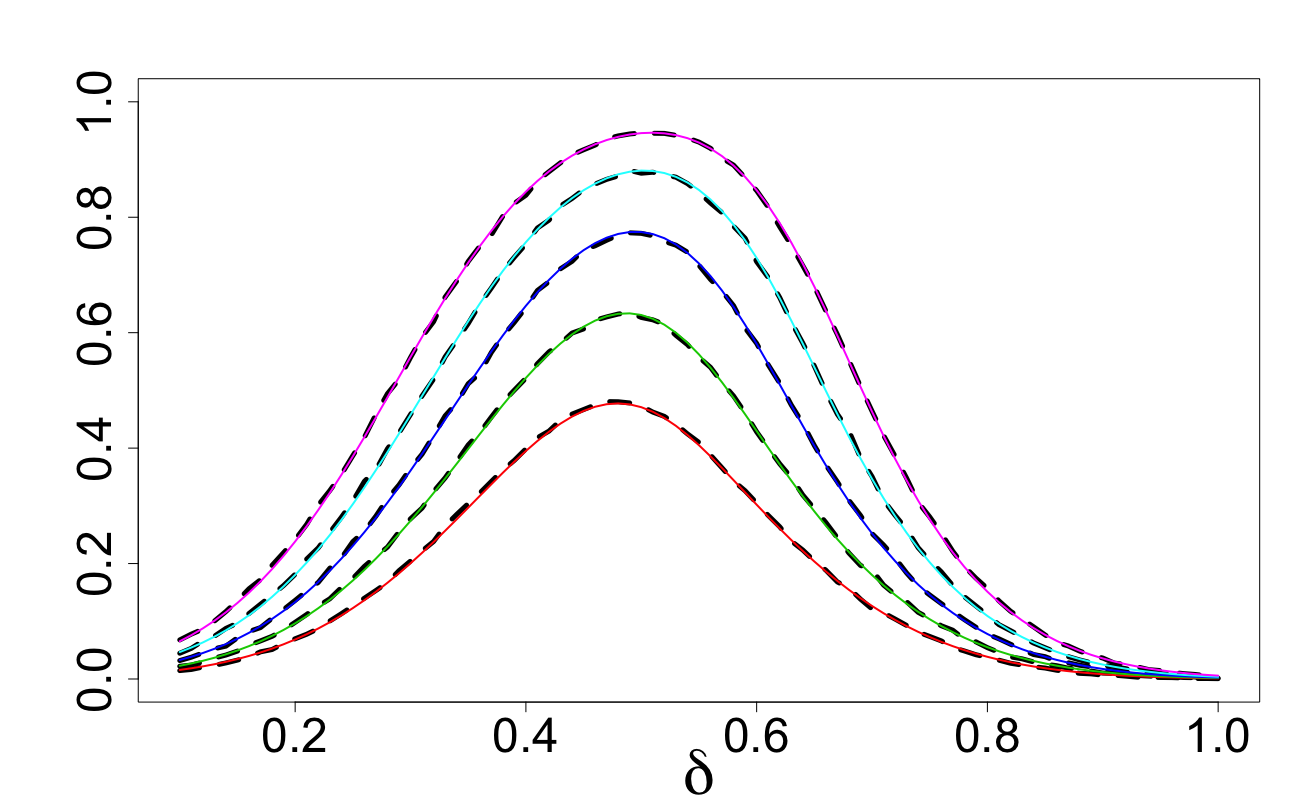}
\caption{$C_d(\mathbb{D}_{n,\delta},r)$ and its approximation: $d=10$, $\\r$ from $0.95$ to $1.15$ increasing by $0.05$}
\label{Covering_approx2}
\end{minipage}
\end{figure}

\begin{figure}[h]
\centering
\begin{minipage}{.5\textwidth}
  \centering
  \includegraphics[width=1\textwidth]{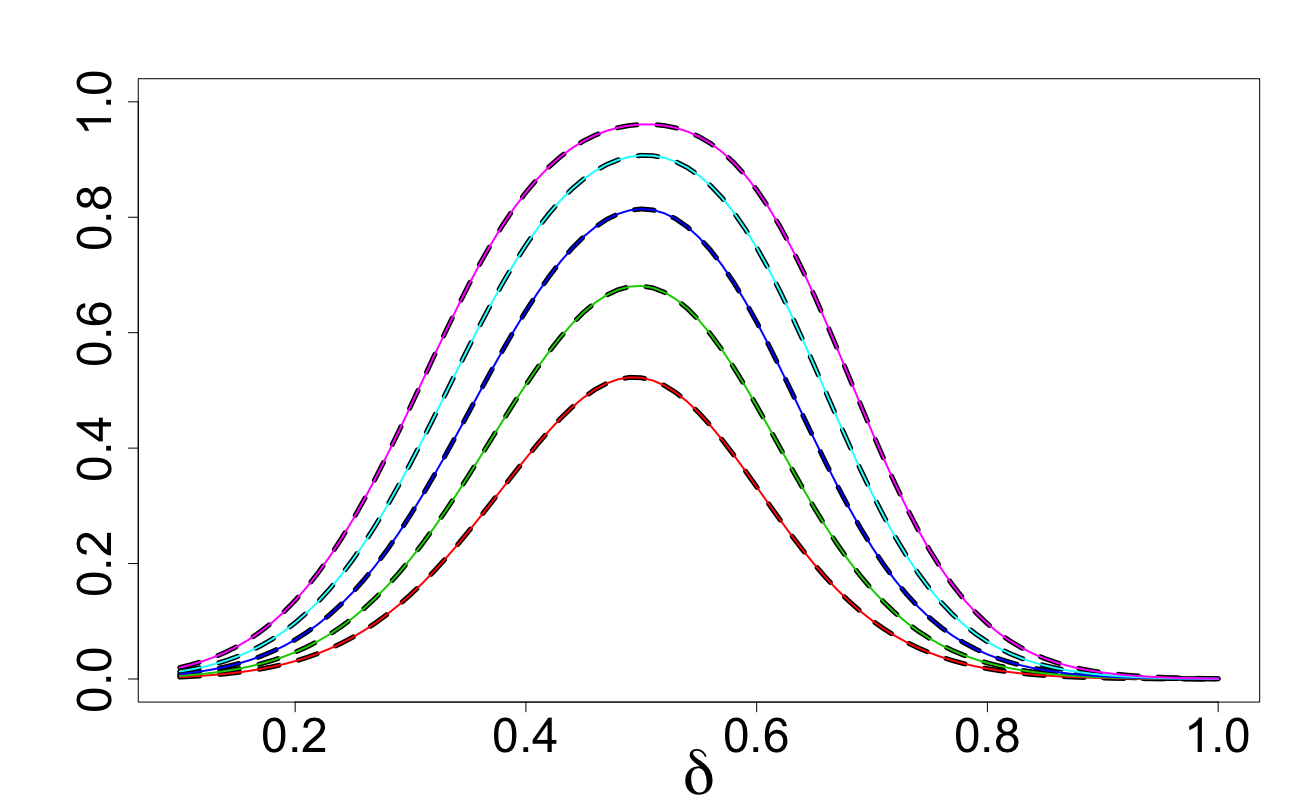}
\caption{$C_d(\mathbb{D}_{n,\delta},r)$ and its approximation: $d=15$, $\\r$ from $1.15$ to $1.35$ increasing by $0.05$ }
\label{Covering_approx5}
\end{minipage}%
\begin{minipage}{.5\textwidth}
  \centering
 \includegraphics[width=1\textwidth]{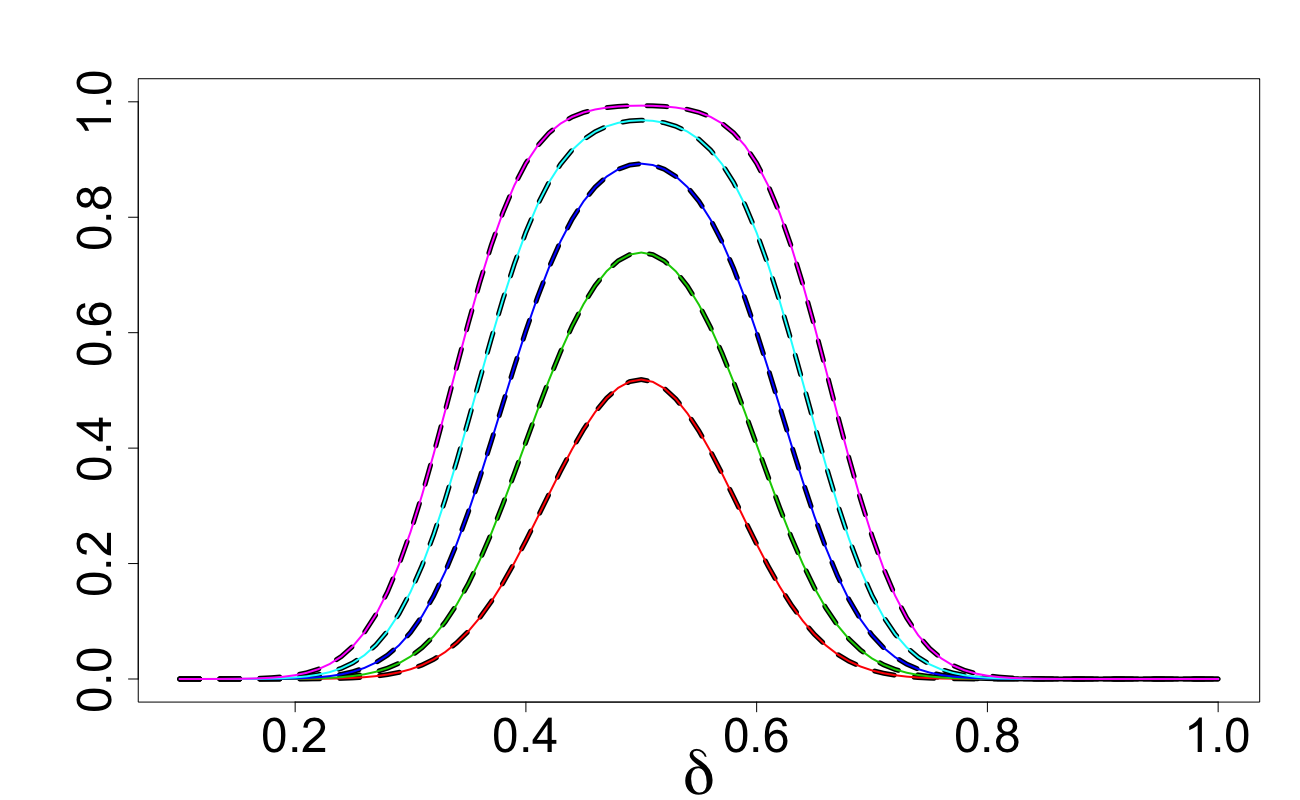}
\caption{$C_d(\mathbb{D}_{n,\delta},r)$ and its approximation:  $d=50$, $\\r$ from $2.05$ to $2.35$ increasing by $0.075$}
\label{Covering_approx6}
\end{minipage}
\end{figure}

In this section, we assess the accuracy of the approximation of $C_d(\mathbb{D}_{n,\delta},r)$ developed in Section~\ref{sec:app1} and the behaviour of $C_d(\mathbb{D}_{n,\delta},r)$ as a function of $\delta$. In Figures~\ref{Covering_approx1} -- \ref{Covering_approx6}, the thick dashed black lines depict $C_d(\mathbb{D}_{n,\delta},r)$ for several different choices of $r$; these values  are obtained via Monte Carlo simulations. The thinner solid lines depict its approximation of Section~\ref{sec:app1}. These figures show that the approximation is extremely accurate for all $r$, $\delta$ and $d$; we emphasise that the approximation remains accurate even for very small dimensions like $d=3$. These figures also clearly demonstrate the $\delta$-effect saying that  a significantly more efficient weak coverage can be achieved with a suitable choice of $\delta$. This is particularly evident in higher dimensions, see Figures~\ref{Covering_approx5} and~\ref{Covering_approx6}.

\begin{figure}[h]
\centering
\begin{minipage}{.5\textwidth}
  \centering
  \includegraphics[width=1\textwidth]{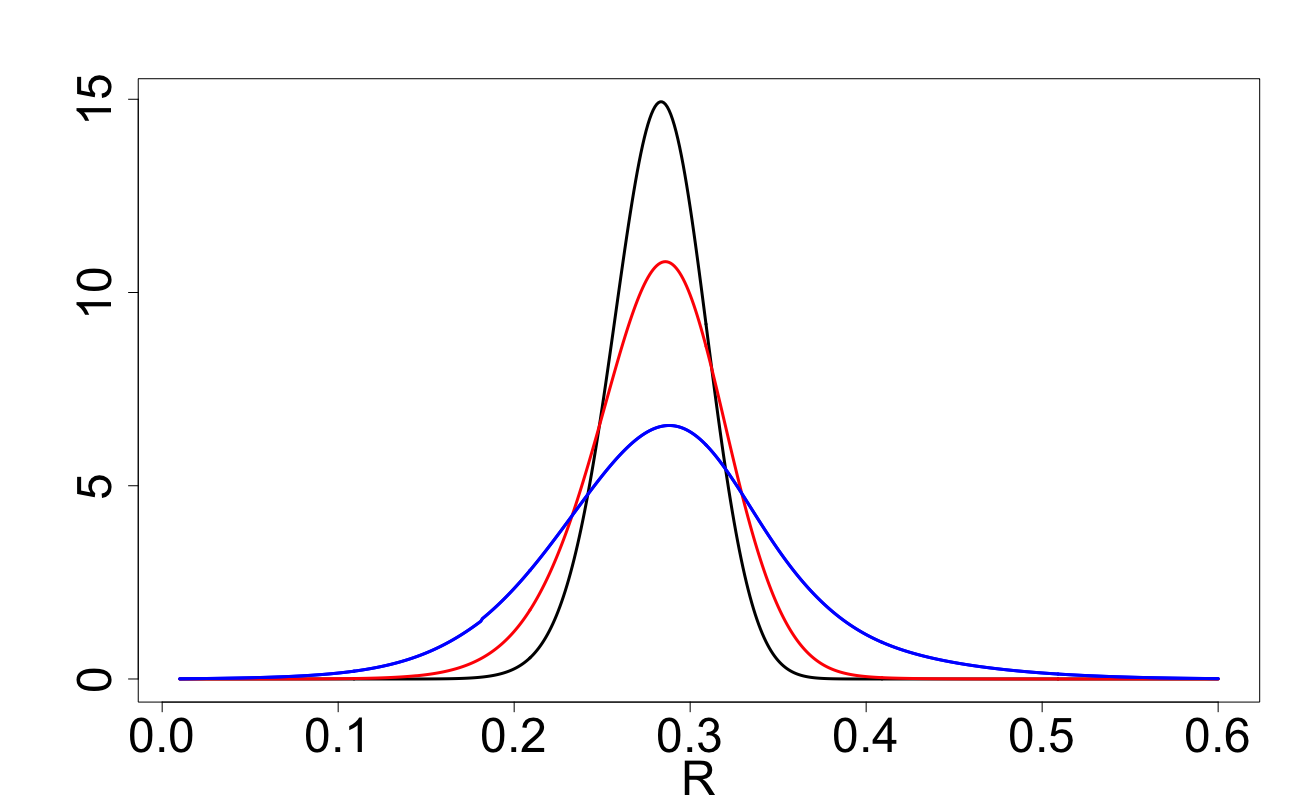}
\caption{Densities $f_d(R)$ for the design $\mathbb{D}_{n,\delta^\ast}$ ; \\ $d=5,10,20$}
\label{density1}
\end{minipage}%
\begin{minipage}{.5\textwidth}
  \centering
 \includegraphics[width=1\textwidth]{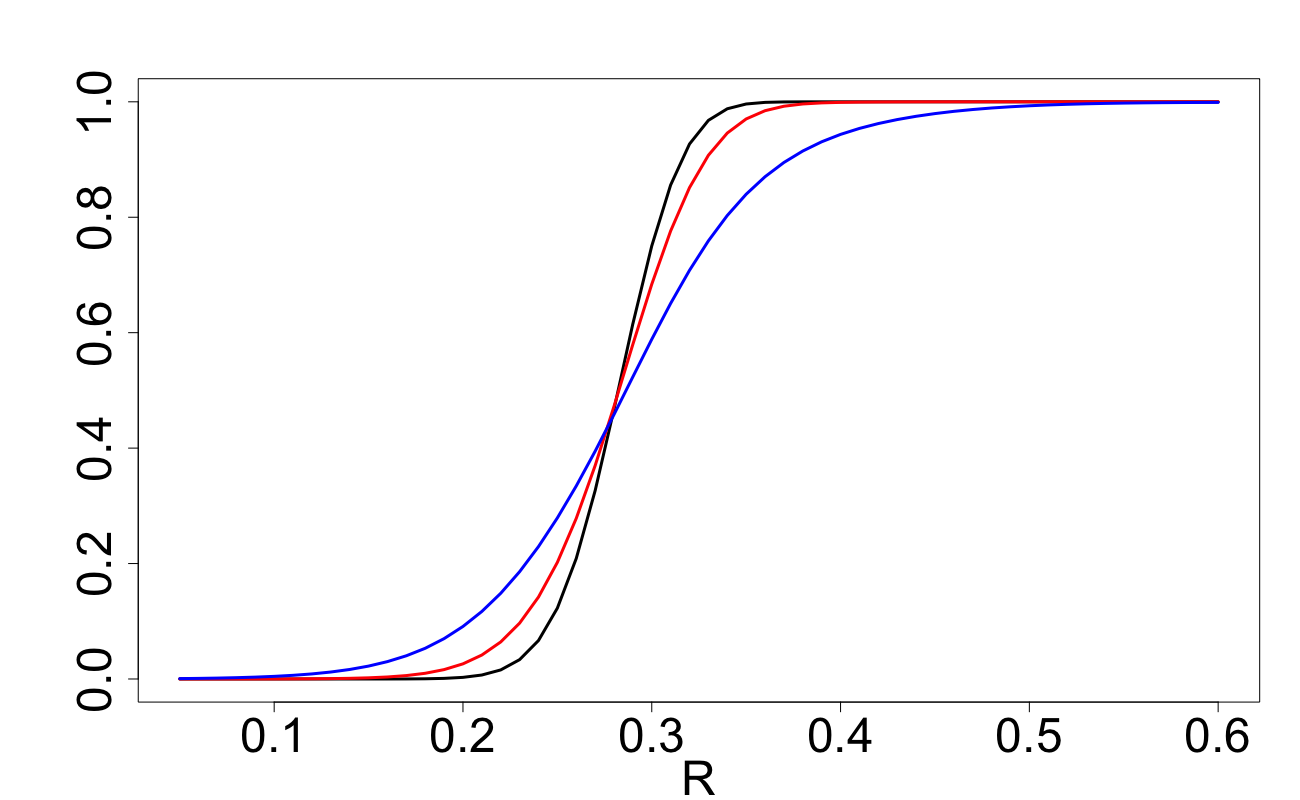}
\caption{c.d.f.'s of $R$ for the design $\mathbb{D}_{n,\delta^\ast}$ ; \\ $d=5,10,20$}
\label{density2}
\end{minipage}
\end{figure}

Figures~\ref{density1} and~\ref{density2} illustrate Theorem~\ref{Asymptotic_theorem} and show the rate of convergence of the covering radii as $d$ increases.
Let the probability density function $f(r)$ be defined by  $d C_d(\mathbb{D}_{n,\delta},r) = f(r)dr$, where $C_d(\mathbb{D}_{n,\delta},r)$ as a function of $r$ is  viewed as the c.d.f. of the r.v. $r=\varrho(X,\mathbb{Z}_n)$, see Section~\ref{sec:relation}.
Trivial calculations yield that the density for the normalised radius $R$ expressed by \eqref{eq:Rr} is  $p_d(R) := {2\sqrt{d}}{n^{-1/d}} f \left( {2\sqrt{d}}{n^{-1/d}}R \right )$. In Figure~\ref{density1}, we depict the density $p_d(\cdot)$ for $d=5,10$ and $20$ with blue, red and black lines respectively. The respective c.d.f.'s $\int_0^R p_d(\tau)d\tau$ are shown in Figure~\ref{density2} under the same colouring scheme.

\subsection{Stochastic dominance}

\begin{figure}[h]
\centering
\begin{minipage}{.5\textwidth}
  \centering
  \includegraphics[width=1\textwidth]{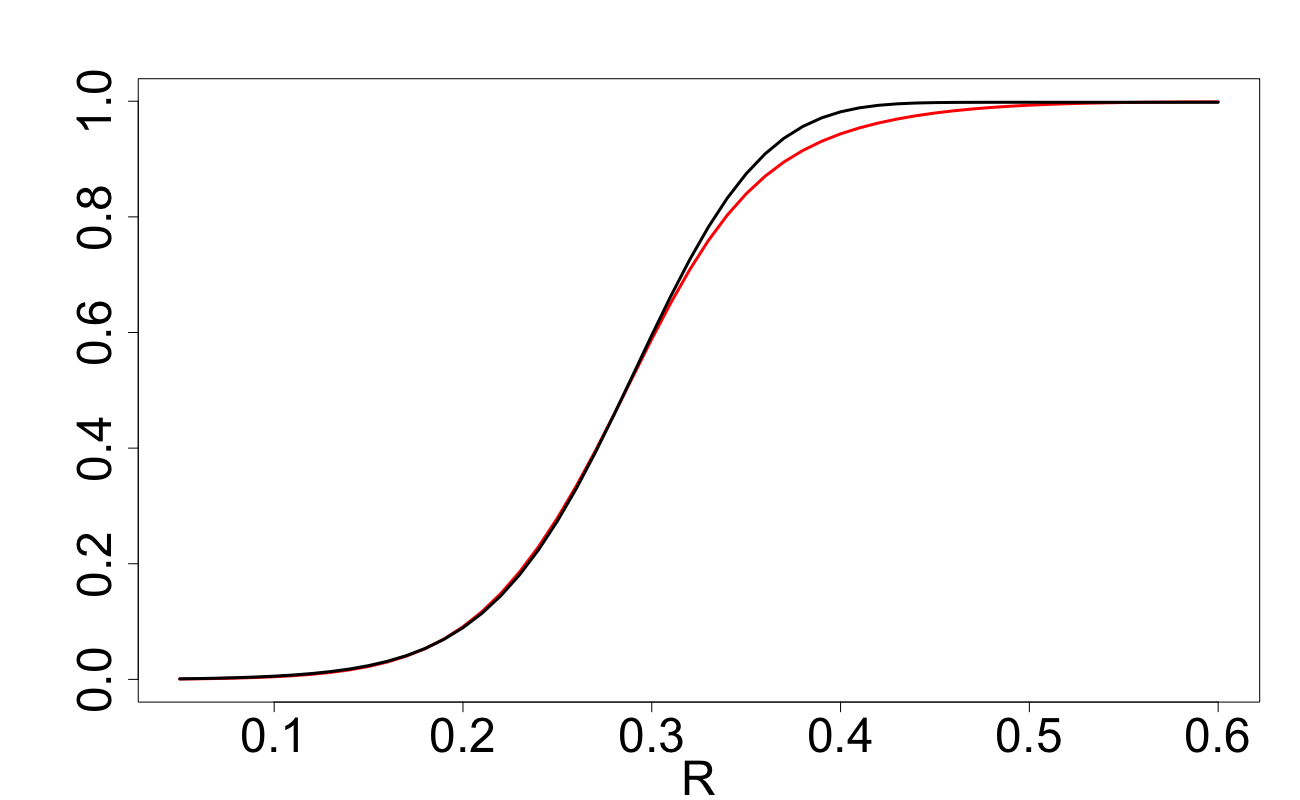}
\caption{$d=5$: design $\mathbb{D}_{n}^{(0)}$ is preferable to \\design  $\mathbb{D}_{n,\delta^\ast}$
}
\label{dom1}
\end{minipage}%
\begin{minipage}{.5\textwidth}
  \centering
 \includegraphics[width=1\textwidth]{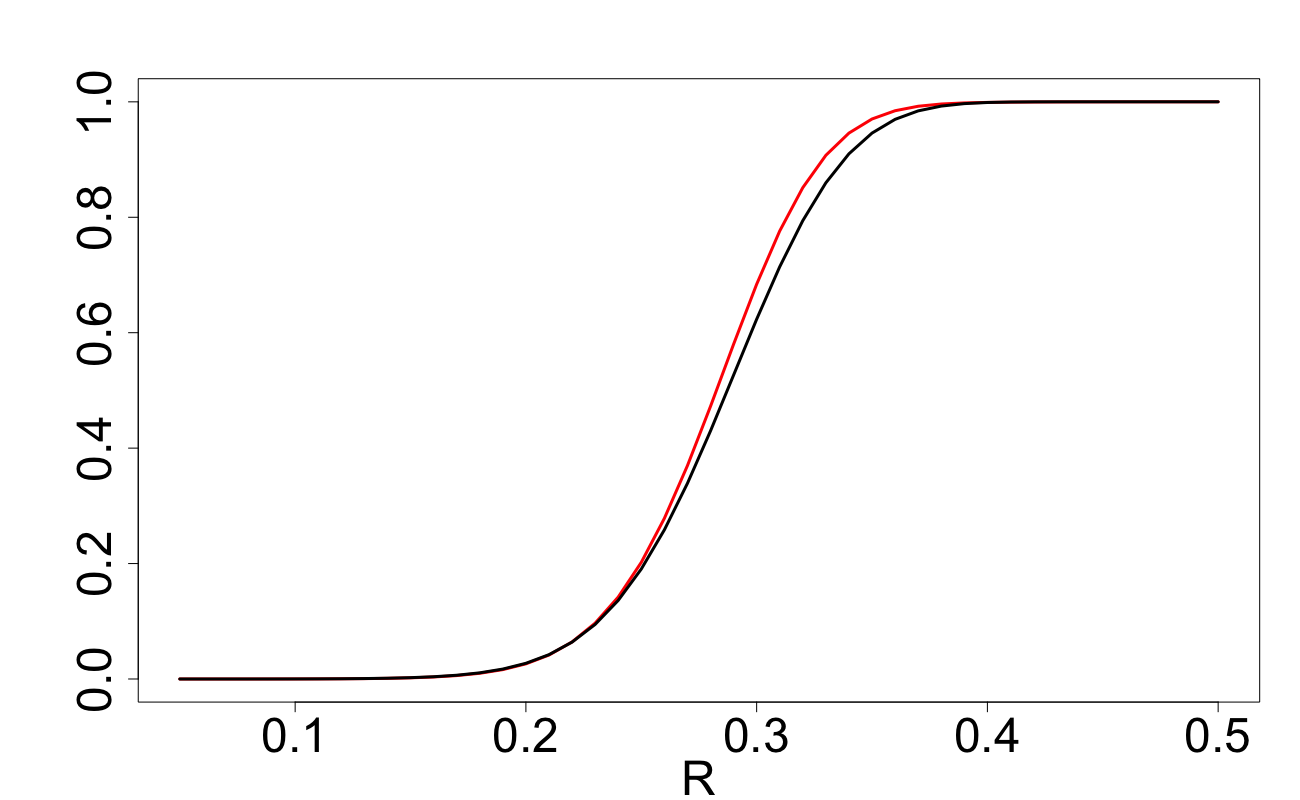}
\caption{$d=10$: design $\mathbb{D}_{n,\delta^\ast}$ stochastically \\dominates design $\mathbb{D}_{n}^{(0)}$
}
\label{dom2}
\end{minipage}
\end{figure}

In Figures~\ref{dom1} and \ref{dom2}, we depict the c.d.f.'s for the normalized distance ${n^{1/d}} \varrho(X,\mathbb{Z}_n)/({ 2 \sqrt{d}} )$  for two designs: $\mathbb{D}_{n,\delta^\ast}$ in red, and $\mathbb{D}_{n}^{(0)}$ in black. We can see that the design $\mathbb{D}_{n,\delta^\ast}$ stochastically dominates  the design $\mathbb{D}_{n}^{(0)}$  for $d=10$ but for $d=5$
the design $\mathbb{D}_{n}^{(0)}$  is preferable to the design $\mathbb{D}_{n,\delta^\ast}$ although there is no clear domination; this is in line with findings from Sections~\ref{quant_exact_results} and \ref{sec:compar}, see e.g. Figure~\ref{Quant_example1},
Tables~\ref{Quant_table} and \ref{Covering_table}.

\begin{figure}[!h]
\centering
\begin{minipage}{.5\textwidth}
  \centering
  \includegraphics[width=1\textwidth]{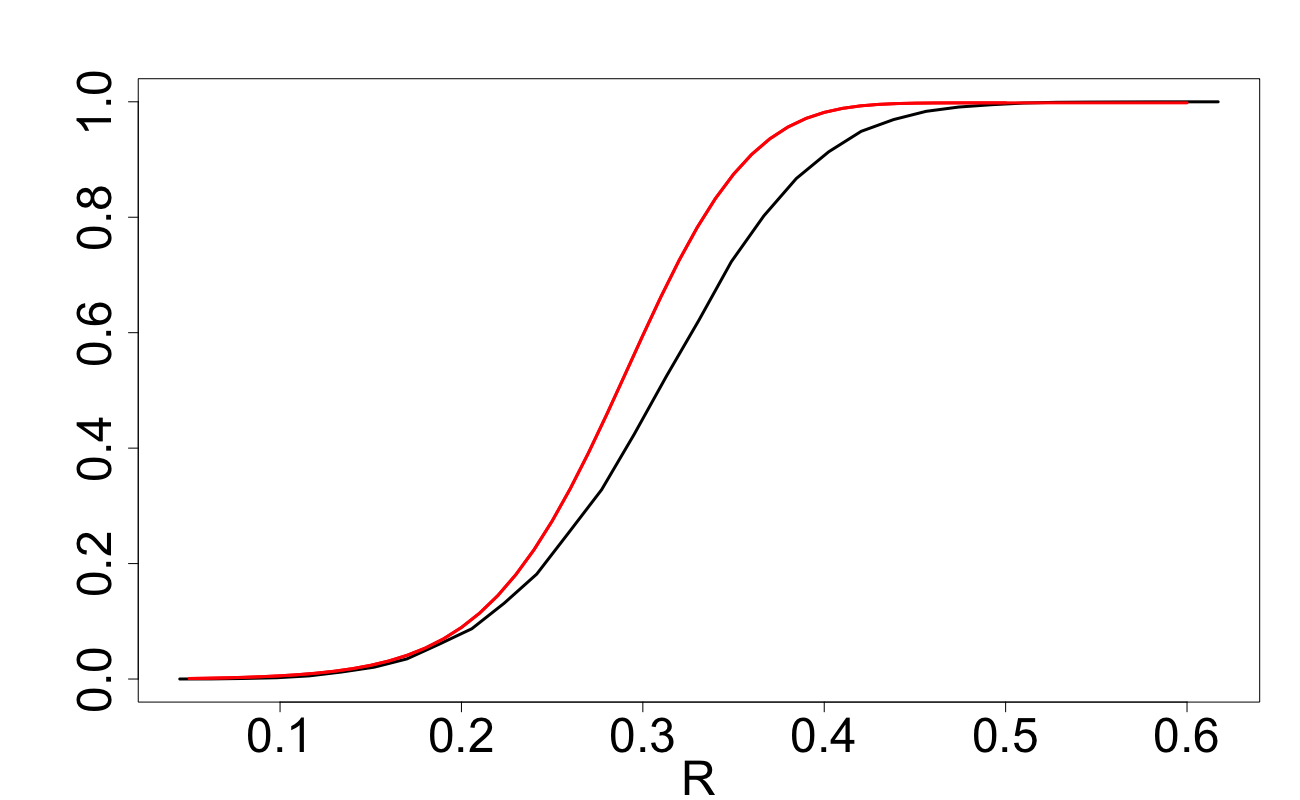}
\caption{$d=5$: design $\mathbb{D}_{n}^{(0)}$ stochastically \\dominates design $\mathbb{S}_{n}$}
\label{dom3}
\end{minipage}%
\begin{minipage}{.5\textwidth}
  \centering
\includegraphics[width=1\textwidth]{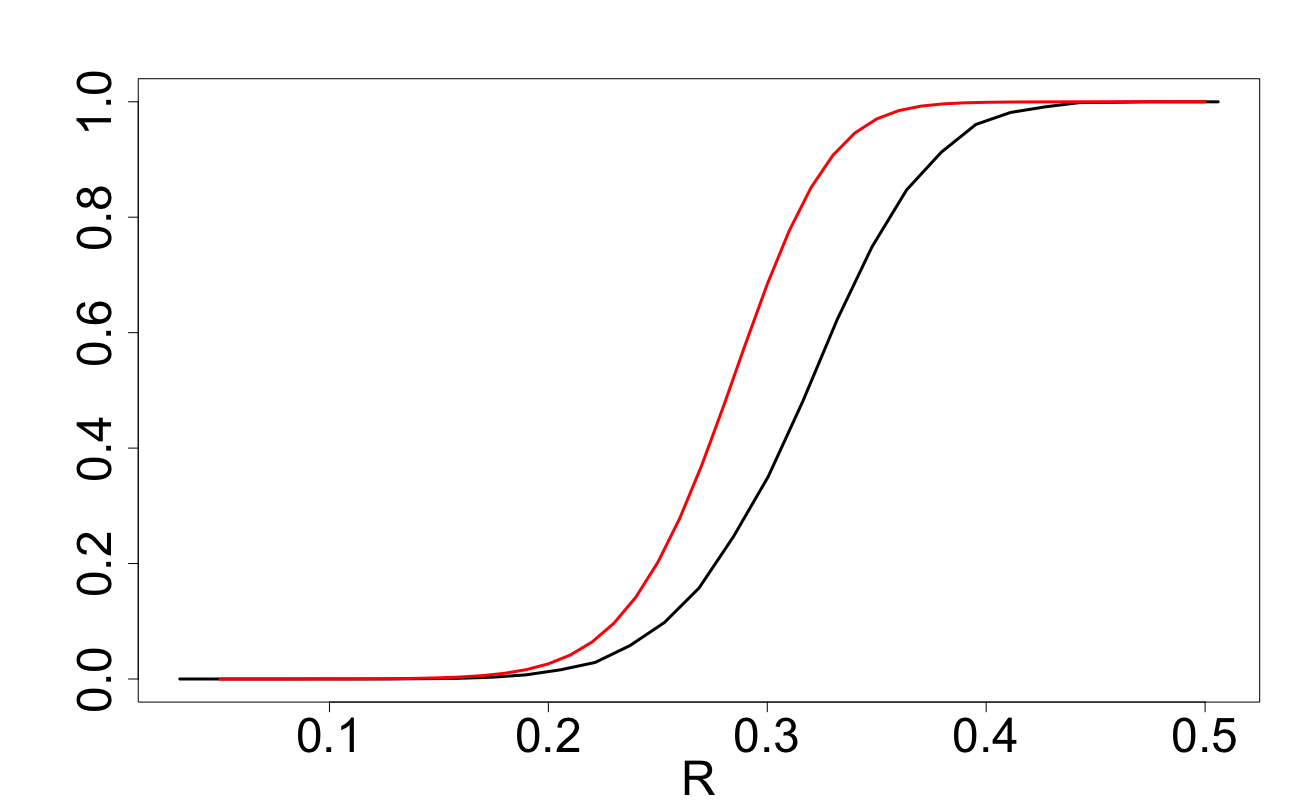}
\caption{$d=10$: design $\mathbb{D}_{n,\delta^\ast}$ stochastically \\dominates design $\mathbb{S}_{n}$}
\label{dom4}
\end{minipage}
\end{figure}

In Figure~\ref{dom3}, we depict the c.d.f.'s for the normalized distance ${n^{1/d}} \varrho(X,\mathbb{Z}_n)/({ 2 \sqrt{d}} )$  for  design $\mathbb{D}_{n}^{(0)}$ (in red) and design $\mathbb{S}_{n}$ (in black). We can see that for $d=5$, the design $\mathbb{D}_{n}^{(0)}$ stochastically dominates  the design $\mathbb{S}_{n}$. The style of Figure~\ref{dom4} is the same as figure Figure~\ref{dom3}, however we set $d=10$ and the design $\mathbb{D}_{n}^{(0)}$ is replaced with the design $\mathbb{D}_{n,\delta^\ast}$. Here we see a very clear stochastic dominance of the design $\mathbb{D}_{n,\delta^\ast}$ over the design $\mathbb{S}_{n}$. All findings are consistent with findings from Section~\ref{sec:compar}, see Tables~\ref{Quant_table} and \ref{Covering_table}.

\section*{Appendix A: Proof of Theorem~\ref{Main_theorem}}\label{Proof_main_thm}


In view of \eqref{voroni_ratio}, $C_d(\mathbb{D}_{n,\delta},r) = V_{d,\boldsymbol{\delta} ,{ r }}$ for all $0 \leq \delta \leq 1$ and $r \geq 0$ and we shall derive expressions for
$ V_{d,\boldsymbol{\delta} ,{ r }}$ rather than $C_d(\mathbb{D}_{n,\delta},r)$.


{$\mathbf{Case (a):}$ $r\leq \delta$}.

\AZ{To prove this case, we observe i) for this range of $r$, ${\cal B}_d(\boldsymbol{\delta},{ r })  \subset [0,1]^d$; ii) the fraction of a cube covered by a ball is
preserved under invertible affine transformations; iii) the affine transformation $x\rightarrow  2x - \boldsymbol{1}$
maps the ball ${\cal B}_d(\boldsymbol{\delta},{ r })$ and cube $[0,1]^d$ to ${\cal B}_d(\boldsymbol{2\delta-1},{2 r })$ and $[-1, 1]^d$, respectively. This
leads to
\bea
V_{d,\boldsymbol{\delta},{ r }}=\frac{{\rm vol}({\cal B}_d(\boldsymbol{\delta},{ r }))}{2\,{\rm vol}([0,1]^d)}= \frac{{\rm vol}({\cal B}_d(\boldsymbol{2\delta-1},{ 2r }))}{2\,{\rm vol}([-1,1]^d)}=\frac12 C_{d,\boldsymbol{2\delta - 1},{2r }} \,.
\eea}
%
%
%
{$\mathbf{Case (b)}$:  $\delta \leq r \leq 1+\delta$}

Using \eqref{eq:uj1} we obtain
\bea
V_{d,\boldsymbol{\delta},{ r }} &=& \frac12 \bigg[ {\rm vol}\left(
{\cal B}_d(\boldsymbol{\delta},{ r }) \cap C_0 \right)         +
d\cdot{\rm vol}\left(
{\cal B}_d(\boldsymbol{\delta},{ r }) \cap U_1 \right) \bigg].
\eea
The first quantity in the brackets has been considered in case (a) and  it is simply $C_{d,\boldsymbol{2\delta - 1},{2r }} $. \AZ{Therefore we aim to reformulate the second quantity within the brackets, $  {\rm vol}\left(
{\cal B}_d(\boldsymbol{\delta},{ r }) \cap U_1 \right)$.
Denote by ${\cal P}(t) = \{ (x_1,x_2,\ldots,x_d):x_1=t \}$, the $(d-1)$-dimensional hyperplane. Then
\bea
{\rm vol}\left( {\cal B}_d(\boldsymbol{\delta},{ r }) \cap U_1 \right) = \int_{\delta-r}^{0}{\rm vol}_{d-1}( {\cal P}(t)\cap {\cal B}_d(\boldsymbol{\delta},{ r }) \cap U_1 )   dt \,.
\eea
Notice further that
\be\label{referee_remark}
 \begin{aligned}
U_1  \cap {\cal P}(t)  = \{ t \} \times [|t|,1]^{d-1}, && \text{ for } -1 \leq t \leq 0 \\
{\cal B}_d(\boldsymbol{\delta},{ r }) \cap {\cal P}(t)  = \{ t \} \times {\cal B}_{d-1}\left(\boldsymbol{\delta},{ \sqrt{r^2-(t-\delta)^2} } \right)&&  \text{ for } \delta-r \leq t \leq 0, r\geq \delta,
  \end{aligned}
\ee
where $\boldsymbol{\delta} =(\delta, \ldots,\delta) \in \mathbb{R}^{d-1}$ and the natural identification of ${\cal P}(t)$ with $\mathbb{R}^{d-1}$ is used.  The r.h.s. in \eqref{referee_remark} are a $(d-1)-$dimensional cube and ball respectively. Since covered fraction is preserved under affine transformations in $\mathbb{R}^{d-1}$, it suffices to construct one, denote by $\phi$, for which $\phi([|t|,1]^{d-1} ) = [-1,1]^d$. In ${\cal P}(t)$, such $\phi$ maps the cube from \eqref{referee_remark} to the standard cube $[-1,1]^d$. Clearly, $\phi$ can be taken as
\bea
\phi: x \rightarrow \frac{x-(\boldsymbol{1+|t|})/2}{({1-|t|})/2} = \frac{2x - (\boldsymbol{1+|t|})}{{1-|t|}} \,,
\eea
where ${\boldsymbol{1}} = (1,\ldots,1)$ and $\boldsymbol{|t|} = (|t|,\ldots,|t|)$ are constant vectors in $\mathbb{R}^{d-1}$. Note that,
\bea
\phi\left({\cal B}_{d-1}\left(\boldsymbol{\delta},{ \sqrt{r^2-(t-\delta)^2} } \right)   \right) = {\cal B}_{d-1}\left( \frac{\boldsymbol{2\delta-(1+|t|)}}{\boldsymbol{1-|t|}},\frac{ 2\sqrt{r^2-(t-\delta)^2} }{1-|t| } \right) \,.
\eea
Finally, by the preservation of covered fraction, we obtain
\bea
{\rm vol}_{d-1}( {\cal P}(t)\cap {\cal B}_d(\boldsymbol{\delta},{ r }) \cap U_1 )= {\rm vol}_{d-1}([|t|,1]^{d-1}) \cdot C_{d-1,\boldsymbol{\frac{2\delta-1- |{t}| }{1-|{t}|}},{ \frac{2\sqrt{r^2-(t-\delta)^2}}{1-|t|} }}  \,.
\eea}
As a result,
\be
V_{d,\boldsymbol{\delta},{ r }} &=& \frac12 \left[ C_{d,\boldsymbol{2\delta - 1},{2r }} + d \int_{\delta-r}^{0} C_{d-1,\boldsymbol{\frac{2\delta-1- |{t}| }{1-|{t}|}},{ \frac{2\sqrt{r^2-(t-\delta)^2}}{1-|t|} }} (1-|t|)^{d-1}\, dt  \right] \label{Case_b_int}\\
&=& \frac12 \left[ C_{d,\boldsymbol{2\delta - 1},{2r }} + d \int_{0}^{r-\delta} C_{d-1,\boldsymbol{\frac{2\delta-1-{t} }{1-{t}}},{ \frac{2\sqrt{r^2-(t+\delta)^2}}{1-t} }} (1-t)^{d-1}\, dt  \right] \nonumber \,.
\ee

 {$\mathbf{Case (c)}$:  $r \geq 1+\delta$:}

Case (c) is almost identical to Case (b), with the only change occurring within the lower limit of integration in \eqref{Case_b_int}; the lower limit of the integral remains at $-1$ for all $r\ge 1+\delta$. Since the steps are almost identical to Case (b), they are omitted and we simply conclude:
\bea
V_{d,\boldsymbol{\delta},{ r }} = \frac12 \left[ C_{d,\boldsymbol{2\delta - 1},{2r }} + d \int_{0}^{1} C_{d-1,\boldsymbol{\frac{2\delta-1- {t} }{1-{t}}},{ \frac{2\sqrt{r^2-(t+\delta)^2}}{1-t} }} (1-t)^{d-1}\, dt  \right] \,.
\eea
 \hfill $\Box$

\section*{Appendix B: Proof of Lemma \ref{Upper_and_lower_bounds}}

(a) Let us first prove the upper bound in \eqref{eq:bounds}. Consider the set $U_j$ defined in \eqref{eq:uj} and the associated set
\bea
U_j^\prime = \left\{  X = (x_1,x_2,\ldots,x_d)\! \in\! [0,1]^d\! :  |x_j|\leq x_k \leq 1\; \mbox{ \rm for all } k\neq j\right\} \subset C_0 \, .
\eea
We have vol($U_j$)=vol($U_j^\prime$)$=1/d$ and
$$
V(\boldsymbol{\delta})= C_0  \bigcup \left[  \bigcup_{j=1}^d U_j\right],\;\;\;\bigcup_{j=1}^d U_j^\prime = C_0=[0,1]^d
$$

Let us prove that for any $r \geq 0$ we have vol($U_j \cap {\cal B}_d(\boldsymbol{\delta},{ r }) $)$\leq$ vol($U_j^\prime \cap {\cal B}_d(\boldsymbol{\delta},{ r }) $).

With any point $X  = (x_1,x_2,\ldots,x_d)\! \in\! U_1^\prime$,
we associate the point $X^\prime  = (-x_1,x_2,\ldots,x_d)\! \in\! U_1$ by simply changing the sign in the first component. For these two points, we have
\bea
\| \boldsymbol{\delta}-X\|^2= \left(x_1-\delta \right)^2+ \sum_{k=2}^{d}\left(x_k-\delta \right)^2< \left(-x_1-\delta \right)^2+ \sum_{k=2}^{d}\left(x_k-\delta \right)^2= \| \boldsymbol{\delta}-X^\prime\|^2
\eea
Therefore, $\| \boldsymbol{\delta}-X\|^2\leq r \Rightarrow \| \boldsymbol{\delta}-X^\prime\|^2\leq r$ yielding:

\be\label{Upper_bound}
{\rm vol}(U_j \cap {\cal B}_d(\boldsymbol{\delta},{ r }) )\leq {\rm vol}(U_j^\prime \cap {\cal B}_d(\boldsymbol{\delta},{ r }) ) \, .
 \ee

To prove the upper bound in \eqref{eq:bounds} for all $r$ we must consider two cases:  $ r\leq \delta$  and $\, r  \geq \delta$.

For $r \leq \delta$, we clearly have
\bea
V_{d,\boldsymbol{\delta},{ r }}  = \frac12 C_{d,\boldsymbol{2\delta - 1},{2r }} \leq  C_{d,\boldsymbol{2\delta - 1},{2r }}
\eea

For $r \geq \delta$,using \eqref{Upper_bound} we have
\bea
V_{d,\boldsymbol{\delta},{ r }} &=& \frac12 \bigg[ {\rm vol}\left(
{\cal B}_d(\boldsymbol{\delta},{ r }) \cap C_0 \right)         +
d\cdot{\rm vol}\left(
{\cal B}_d(\boldsymbol{\delta},{ r }) \cap U_1 \right) \bigg] \\
&\leq& \frac12 \bigg[ {\rm vol}\left(
{\cal B}_d(\boldsymbol{\delta},{ r }) \cap C_0 \right)         +
d\cdot{\rm vol}\left({\cal B}_d(\boldsymbol{\delta},{ r }) \cap U'_1 \right) \bigg]  \\
&=&   {\rm vol}({\cal B}_d(\boldsymbol{\delta},{ r }) \cap C_0 )\\
&=&   C_{d,\boldsymbol{2\delta - 1},{2r }}
\eea

 and hence the upper bound in \eqref{eq:bounds}.

(b) Consider now the lower bound in  \eqref{eq:bounds}. For $j\geq2$, with the set $U_j$ we now associate the set
\bea
V_j = \left\{ \tilde{ X }= (x_1,\ldots,x_d)\! : -1\leq x_1 \leq 0 ,\; 0\leq x_m \leq 1\; (\mbox{ \rm for } m>1), \;|x_j|\leq |x_k| \leq 1\; \mbox{ \rm for } k\neq j\right\} \, .
\eea

 With any point $X  = (x_1,x_2,\ldots,x_d)\! \in\! U_j$ (here $x_j$ is negative and $|x_j|\leq |x_k| \leq 1\; \mbox{ \rm for } k\neq j $)
 we associate point $\tilde{X}  = (-x_1,x_2,\ldots,x_{j-1}, -x_j, x_{j+1}, \ldots, x_d)\! \in\! V_j$ by changing sign in the first and $j$-the component of $X \in U_j$.

Setting without loss of generality $j=2$, we have for these two points:
\bea
\| \boldsymbol{\delta}-X\|^2&=& \left(x_1-\delta \right)^2+ \left(x_2-\delta \right)^2+  \sum_{k=3}^{d}\left(x_k-\delta \right)^2\\
& \leq & \left(-x_1-\delta \right)^2+ \left(-x_2-\delta \right)^2+  \sum_{k=3}^{d}\left(x_k-\delta \right)^2= \| \boldsymbol{\delta}-\tilde{X}\|^2\, ,
\eea
where the inequality follows from the inequalities $x_1\geq 0$, $x_2<0$ and $|x_2|<x_1$ containing in the definition of $U_2$.

Therefore, $\| \boldsymbol{\delta}-\tilde{X}\|^2\leq r \Rightarrow \| \boldsymbol{\delta}-{X}\|^2\leq r$ implying:
\be\label{Lower_bound}
{\rm vol}(U_j \cap {\cal B}_d(\boldsymbol{\delta},{ r }) )\geq {\rm vol}(V_j \cap {\cal B}_d(\boldsymbol{\delta},{ r }) )\, .
 \ee

To prove the lower bound in \eqref{eq:bounds} for all $r$ we must consider two cases: $ r\leq \delta$  and $ \, r  \geq \delta$.

For $r \geq \delta$, using \eqref{Lower_bound} we  have
\bea
V_{d,\boldsymbol{\delta},{ r }} &=& \frac12 \bigg[ {\rm vol}\left(
{\cal B}_d(\boldsymbol{\delta},{ r }) \cap C_0 \right)         +
\sum_{i=1}^{d}{\rm vol}\left(
({\cal B}_d(\boldsymbol{\delta},{ r }) \cap U_i )\right) \bigg] \\
&\geq& \frac12 \bigg[ {\rm vol}\left(
{\cal B}_d(\boldsymbol{\delta},{ r }) \cap C_0 \right)         +
{\rm vol}({\cal B}_d(\boldsymbol{\delta},{ r }) \cap U_1 )+ \sum_{i=2}^{d}
\left({\cal B}_d(\boldsymbol{\delta},{ r }) \cap V_i )\right) \bigg] \\
&=& \frac12 \bigg[ {\rm vol}\left(
{\cal B}_d(\boldsymbol{\delta},{ r }) \cap C_0 \right)         +
{\rm vol}({\cal B}_d(\boldsymbol{\delta},{ r }) \cap C_1)  \bigg]  \, ,\\
\eea
where $C_1$ is given in \eqref{adjacent_cubes}.
\AZ{To compute ${\rm vol}({\cal B}_d(\boldsymbol{\delta},{ r }) \cap C_1)$, we shall use a similar technique to the proof of Theorem~\ref{Main_theorem}. The affine transformation
\bea
x \rightarrow 2x +  (1,-1,-1,\ldots,-1)
\eea
maps the ball ${\cal B}_d(\boldsymbol{\delta},{ r })$ and the cube $C_1$ to ${\cal B}_d({A},{2r })$ and $[-1,1]^d$ respectively, where ${A}= \left(2\delta+1,2\delta-1,\ldots,2\delta-1 \right)$. Since the fraction of covered volume is preserved under invertible affine transformations, one has
%
\bea
\frac{{\rm vol}({\cal B}_d(\boldsymbol{\delta},{ r }) \cap C_1)}{{\rm vol}(C_1)} = C_{d,A,2r}
\eea}
and hence we can conclude:
\bea
V_{d,\boldsymbol{\delta},{ r }} &\geq& \frac12 \bigg[ {\rm vol}\left(
{\cal B}_d(\boldsymbol{\delta},{ r }) \cap C_0 \right)         +
{\rm vol}({\cal B}_d(\boldsymbol{\delta},{ r }) \cap C_1)  \bigg] \\
&=&  \frac{1}{2} \left [ C_{d,\boldsymbol{2\delta - 1},{2r }} + C_{d,A,2r}\right] \,.
\eea

For $r \leq \delta$, since ${\rm vol}({\cal B}_d(\boldsymbol{\delta},{ r }) \cap C_1)=C_{d,A,2r} =0$, we  have
\bea
V_{d,\boldsymbol{\delta},{ r }} &=&  \frac{1}{2} \left [ C_{d,\boldsymbol{2\delta - 1},{2r }} + C_{d,A,2r}\right] \,\eea

 and hence the lower bound in \eqref{eq:bounds}.
 \hfill $\Box$\\

\section*{Appendix C: Proof of Theorem~\ref{Asymptotic_theorem}.}

Before proving Theorem~\ref{Asymptotic_theorem}, we prove three auxiliary lemmas.



\begin{lemma}\label{lim_of_c} Let $r=r_{\alpha,d}=\alpha \sqrt{d}$ with $\alpha \geq 0$ and $Z_{a,b;d}=(a,b,b ,\ldots, b)\in \mathbb{R}^d$. 
Then the limit  $ \lim_{d \to \infty } C_{d,{Z_{a,b;d}},{2r }}$ exists and
\bea
\lim_{d \to \infty} C_{d,Z_{a,b;d},{2r }} = \begin{cases}
0 \text{ if    } \alpha<\frac{1}{2}\sqrt{\frac13+b^2}\\
1/2  \text{ if    } \alpha=\frac{1}{2}\sqrt{\frac13+b^2}\\
1 \text{ if    } \alpha>\frac{1}{2}\sqrt{\frac13+b^2}\,
\end{cases}
\eea
\end{lemma}

{ \bf Proof.}
Define
\bea
t_\alpha =  \frac{\sqrt{3}(d(4\alpha^2 -b^2-1/3) +b^2-a^2)}{2\sqrt{a^2+(d-1)b^2+d/15}}\, .
\eea


As the r.v. $\eta_z$ introduced in Appendix D are concentrated on a finite interval, for finite $a$ and $b$  the quantities of $\rho_{a} :=  \mathbb{E}(|\eta_a - a^2 -\frac1{3}|^3)$ and $\rho_{b} :=  \mathbb{E}(|\eta_b - b^2 -\frac1{3}|^3)$  are  bounded.
By applying Berry-Esseen theorem (see  \S 2, Chapter 5 in \cite{petrov2012sums}) to $C_{d,Z_{a,b},{2r }}$, there exists some constant $C$ such that
\bea
-\frac{C \cdot \max\{ \rho_{a}/\sigma_{a}^2,\rho_{b}/\sigma_{b}^2 \} }{\left (\sigma_{a}^2+(d-1)\sigma_{b}^2  \right)^{1/2}}+\Phi \left( t_\alpha \right)\leq C_{d,Z_{a,b},{2r }}\leq \Phi \left( t_\alpha \right)+\frac{C \cdot \max\{ \rho_{a}/\sigma_{a}^2,\rho_{b}/\sigma_{b}^2 \} }{\left (\sigma_{a}^2+(d-1)\sigma_{b}^2  \right)^{1/2}} \, ,
\eea
where $\sigma_{a}^2={\rm var}(\eta_a)$ and  $\sigma_{b}^2={\rm var}(\eta_b)$.
By the squeeze theorem, it is clear that if $4\alpha^2 -b^2-1/3>0$ and hence $\alpha > \frac{1}{2}\sqrt{\frac13+b^2}$, then $ C_{d,Z_{a,b},{2r }}\rightarrow 1$ as $d \rightarrow \infty$. If $\alpha < \frac{1}{2}\sqrt{\frac13+b^2}$, then $ C_{d,Z_{a,b},{2r }}\rightarrow 0$ as $d \rightarrow \infty$. If $\alpha =\frac{1}{2}\sqrt{\frac13+b^2}$, then $ C_{d,Z_{a,b},{2r }}\rightarrow 1/2$ as $d \rightarrow \infty$.
\hfill $\Box$\\

\begin{lemma}\label{Opt_alpha}
Let $r = \alpha\sqrt{d}$. Then for $\boldsymbol{\delta}= (\delta,\delta,\ldots,\delta)$, we have:
\bea
\lim_{d \to \infty} V_{d,\boldsymbol{\delta},{ r }} = \lim_{d \to \infty} C_{d,\boldsymbol{2\delta-1},{2r }} = \begin{cases}
0 \text{ if    } \alpha<\frac{\sqrt{1/3 + (2\delta -1)^2}}{2}\\
1/2  \text{ if    } \alpha=\frac{\sqrt{1/3 + (2\delta -1)^2}}{2}\\
1 \text{ if    } \alpha>\frac{\sqrt{1/3 + (2\delta -1)^2}}{2}\,
\end{cases}
\eea

\end{lemma}

{ \bf Proof.}
Using Lemma~\ref{lim_of_c} with $Z_{a,b} =A =\left(2\delta+1,2\delta-1,\ldots,2\delta-1 \right) $, we obtain:
\bea
\lim_{d \to \infty}C_{d,A,{2r }} =\lim_{d \to \infty} C_{d,\boldsymbol{2\delta-1},{2r }} =\begin{cases}
0 \text{ if    } \alpha<\frac{\sqrt{1/3 + (2\delta -1)^2}}{2}\\
1/2  \text{ if    } \alpha=\frac{\sqrt{1/3 + (2\delta -1)^2}}{2}\\
1 \text{ if    } \alpha>\frac{\sqrt{1/3 + (2\delta -1)^2}}{2}\,
\end{cases}
\eea
By then applying the squeeze theorem to the bounds in Lemma~\ref{Upper_and_lower_bounds} using the fact from Lemma~\ref{lem:equal} we have $V_{d,\boldsymbol{\delta},{ r }} =C_d(\mathbb{Z}_n,r) $, we obtain the result.
\hfill $\Box$\\

To determine the value of $r$ that leads to the full coverage, we utilise the following simple lemma.
\AZ{\begin{lemma}\label{r_max}
For design ${\mathbb{D}_{n,\delta}}$, the smallest value of $r$ that ensures a complete coverage of  $[-1,1]^d$ satisfies:
\bea
\lim_{d \rightarrow \infty }\frac{r_{1}}{\sqrt{d}} = \begin{cases}
{(1-\delta)}\text{ if    } \delta \le 1/2\\
{\delta} \text{ if    }  \delta >1/2\end{cases}\,.
\eea
\end{lemma}}

{\bf Proof of Theorem~\ref{Asymptotic_theorem}.}

 From Lemma~\ref{Opt_alpha}, it is clear that the smallest $\alpha$ and hence $r$ is attained with $\delta=1/2$. Moreover, Lemma~\ref{Opt_alpha} provides:
\bea
\lim_{d \to \infty} V_{d,\boldsymbol{1/2},{ r }} = \lim_{d \to \infty} C_{d,\boldsymbol{0},{2r }} = \begin{cases}
0 \text{ if    } \alpha<\frac{1}{2\sqrt{3 }}\\
1/2  \text{ if    } \alpha=\frac{1}{2\sqrt{3 }}\\
1 \text{ if    } \alpha>\frac{1}{2\sqrt{3 }}\,
\end{cases}
\eea
meaning for any $0<\gamma<1$, $r_{1-\gamma} =\frac{\sqrt{d}}{2\sqrt{3 }}$. By then applying Lemma~\ref{r_max} with $\delta=1/2$, we obtain $r_{1} = \sqrt{d}/2$ and hence $r_{1-\gamma}/r_{1} \rightarrow 1/\sqrt{3}$ as $d\rightarrow \infty.$
\hfill $\Box$\\

\section*{Appendix D: Derivation of approximation  \eqref{eq:inters2f_corrected}}
\label{Approx_formulation}

Let $U=(u_1, \ldots, u_d) $ be a random vector with uniform distribution on $[-1,1]^d$ so that $u_1, \ldots, u_d$ are i.i.d.r.v. uniformly distributed on $[-1,1]$. Then for given $Z=(z_1, \ldots, z_d)  \in \mathbb{R}^d$ and any $r>0$,
\bea
C_{d,Z,{ r }}\!= \! \mathbb{P} \left\{ \| U\!-\!Z \|\! \leq \! { r } \right\}\!= \! \mathbb{P} \left\{ \| U\!-\!Z \|^2 \leq  { r^2 } \right\}\!= \! \mathbb{P} \left\{\sum_{j=1}^d (u_j\!-\!z_j)^2 \leq  { r }^2 \right\}   .\;\;
\eea
That is,  $C_{d,Z,{ r }}$, as a function of ${ r }$,   is the c.d.f. of the r.v. $\| U-Z \| $.

Let  $u$ have the  uniform distribution on $[-1,1]$ and $z \in \mathbb{R} $.
The first three central moments  of the r.v. $\eta_z = (u-z)^2$ can be easily computed:
\be
\label{eq:intersf}
\mathbb{E}\eta_z =z^2 +\frac1{3}, \;\;{\rm var}(\eta_z) = \frac43 \left(z^2 +\frac1{15} \right)\, , \;\;
\mu_{z}^{(3)}=E \left[\eta_{z} - E\eta_{z}\right]^3 = \frac{16}{15} \left(z^2 +\frac{1}{63} \right)  \, .
\ee


Consider the r.v.
$
\| U-Z \|^2 =\sum_{i=1}^d \eta_{z_j}=\sum_{j=1}^d (u_j-z_j)^2\, .
$
From \eqref{eq:intersf} and independence of $u_1, \ldots, u_d$, we obtain
\bea
\mu_{d,Z}=\mathbb{E}\| U-Z \|^2  =\|Z\|^2 +\frac{d}{3}\, \;\;
 {\sigma}_{d,Z}^2={\rm var}(\| U-Z \|^2 ) = \frac43 \left(\|Z\|^2 +\frac{d}{15}\right)\, \;\;\;\;\;\;
\eea
and
\bea
 {\mu}_{d,Z}^{(3)}= \mathbb{E}\left[\| U-Z \|^2- \mu_{d,Z}\right]^3  = \sum_{j=1}^d     \mu_{z_j}^{(3)} =
 \frac{16}{15} \left(\|Z\|^2 +\frac{d}{63}\right)\, .\;\;\;\;\;\;
\eea

If $d$ is large enough then the conditions of the CLT for $\| U-Z \|^2$ are approximately met and  the distribution of $\| U-Z \|^2 $
 is approximately normal with mean $\mu_{d,Z}$ and variance ${\sigma}_{d,Z}^2$. That is, we can approximate
$C_{d,Z,{ r }}$
by
\be
\label{eq:inters2f}
C_{d,Z,{ r }} \cong \Phi \left(\frac{{ r }^2-\mu_{d,Z}}{{\sigma}_{d,Z}} \right) \, ,
\ee
where $\Phi (\cdot)$   is the c.d.f. of the standard normal distribution:
$$
\Phi (t) = \int_{-\infty}^t \varphi(v)dv\;\;{\rm with}\;\; \varphi(v)=\frac{1}{\sqrt{2\pi}} e^{-v^2/2}\, .
$$
The approximation \eqref{eq:inters2f} can be improved by
using an Edgeworth-type  expansion in the CLT for sums of independent  non-identically distributed r.v.

General expansion in the central limit theorem for sums of independent non-identical r.v. has been derived
by V.Petrov, see
Theorem 7 in Chapter 6 in  \cite{petrov2012sums}; the first three terms of this expansion have been specialized  in
Section 5.6 in \cite{petrov}.
By using only the first term in this expansion,
we obtain the following approximation for the distribution function of $\| U-Z \|^2 $:
\bea
P\left(\frac{\| U-Z \|^2-\mu_{d,Z}}{\sigma_{d,Z}} \leq x \right) \cong \Phi(x) + \frac{ {\mu}_{d,Z}^{(3)}}{6  ({\sigma}_{d,Z}^2)^{3/2} }(1-x^2)\varphi(x),
\eea
leading to the following improved  form of \eqref{eq:inters2f}:
\bea
C_{d,Z,{ r }} \cong \Phi(t) + \frac{  \|Z\|^2+d/63}{5\sqrt{3} (\|Z\|^2+d/15)^{3/2} }(1-t^2)\varphi(t) \, ,
\eea
where
\bea
t = t_{d,\|Z\|,{ r }}= \frac{{ r }^2-\mu_{d,Z}}{{\sigma}_{d,Z}}
=  \frac{\sqrt{3}(r^2- \|Z\|^2 -d/3)}{2\sqrt{ \|Z\|^2 +{d}/{15}} }\,.\\\nonumber
\eea


\begin{acknowledgements}
The authors are grateful to our colleague Iskander Aliev for numerous intelligent discussions. The authors are very grateful to the referee for their comments that greatly improved the presentation of the paper. Furthermore, we are grateful for the concise proof of Theorem~\ref{Proof_main_thm} suggested by the reviewer which has been included in this manuscript.
\end{acknowledgements}

\bibliographystyle{unsrt}
\bibliography{large_dimension}

\end{document}